\documentclass[journal]{IEEEtran}      

\pdfoutput=1
\IEEEoverridecommandlockouts                              

\usepackage{graphicx}
\usepackage{amsmath,bm}
\usepackage{amssymb}
\usepackage{algorithm}
\usepackage{algpseudocode}
\usepackage{cite}
\usepackage{hyperref}
\usepackage[caption=false]{subfig}
\usepackage{mathtools}
\hypersetup{hidelinks}

\newtheorem{Remark}{Remark}
\newtheorem{Corollary}{Corollary}
\newtheorem{Definition}{Definition}

\newenvironment{Proof}{\noindent{\em Proof:\/}}{\hfill \IEEEQED\par}
\newtheorem{Theorem}{Theorem}
\newtheorem{Lemma}{Lemma}

\newtheorem{Assumption}{Assumption}

\algdef{SE}[DOWHILE]{Do}{doWhile}{\algorithmicdo}[1]{\algorithmicwhile\ #1}%

\newcommand{\mathactivatecomma}{%
  \begingroup\lccode`~=`\,
  \lowercase{\endgroup\edef~}{\mathchar\the\mathcode`\,\penalty0 }}

\algnewcommand{\Initialize}[1]{%
  \State \textbf{Initialize: $j\in \mathcal{V}^i, i \in \mathcal{N}$}
  \Statex \hspace*{\algorithmicindent}\parbox[t]{.8\linewidth}{\raggedright #1}
}

\algnewcommand{\Iteration}[1]{%
  \State \textbf{Iteration $(t\geq 0)$: $j\in \mathcal{V}^i, i \in \mathcal{N}$}
  \Statex \hspace*{\algorithmicindent}\parbox[t]{.8\linewidth}{\raggedright #1}
}

\algnewcommand{\Output}[1]{%
  \State \textbf{Output: $j\in \mathcal{V}^i, i \in \mathcal{N}$}
  \Statex \hspace*{\algorithmicindent}\parbox[t]{.8\linewidth}{\raggedright #1}
}

\title{\LARGE \bf
Distributed Nash Equilibrium Seeking in N-Cluster Games with Fully Uncoordinated Constant Step-Sizes
}

\author{Yipeng Pang and Guoqiang Hu
\thanks{This research was supported by Singapore Ministry of Education Academic Research Fund Tier 1 RG180/17(2017-T1-002-158).}
\thanks{Y. Pang and G. Hu are with the School of Electrical and Electronic Engineering, Nanyang
Technological University, 639798, Singapore
        {\tt\small ypang005@e.ntu.edu.sg, gqhu@ntu.edu.sg}.}%
}

\begin{document}

\bstctlcite{IEEEexample:BSTcontrol}

\maketitle
\thispagestyle{empty}
\pagestyle{empty}

\begin{abstract}
Distributed optimization and Nash equilibrium (NE) seeking problems have drawn much attention in the control community recently. This paper studies a class of non-cooperative games, known as $N$-cluster game, which subsumes both cooperative and non-cooperative nature among multiple agents in the two problems: solving distributed optimization problem within the cluster, while playing a non-cooperative game across the clusters. Moreover, we consider a \textit{partial-decision information} game setup, \textit{i.e.}, the agents do not have direct access to other agents' decisions, and hence need to communicate with each other through a directed graph. To solve the $N$-cluster game problem, we propose a distributed NE seeking algorithm by a synthesis of consensus and gradient tracking. Unlike other existing discrete-time methods for $N$-cluster games where either a common step-size is publicly known by all agents or only known by agents from the same cluster, the proposed algorithm can work with fully uncoordinated constant step-sizes, which allows the agents (both within and across the clusters) to choose their own preferred step-sizes. We prove that all agents' decisions converge linearly to their corresponding NE so long as the largest step-size and the heterogeneity of the step-sizes are small. We verify the derived results through a numerical example in a Cournot competition game.
\end{abstract}

\section{Introduction}
Simultaneous social cost minimization and Nash equilibrium (NE) seeking among multiple clusters (or coalitions) modeled by $N$-cluster game have received great attention in recent researches, due to its wide application in many fields, such as business management, transportation systems, political science, sports, to list a few. In such $N$-cluster games, the agents in the same cluster cooperatively minimize a cluster-level cost function, and collectively act as a virtual player to play an $N$-player non-cooperative game across clusters. If each agent is assumed to have direct access to all other agents' decisions, either by observation or via a virtual central coordinator, it is known as \textit{full-decision information} setup. This problem setup could be restrictive, especially in the scenarios where there is no central node with bidirectional communication with all agents. Hence, a \textit{partial-decision information} setup is motivated, where agents are connected by an underlying communication graph, and they may only access the decisions of their neighbors. 

\textit{Related work:} 
Distributed NE seeking under partial-decision information over graphs have been researched, see \cite{Salehisadaghiani2016,DePersis2019,Pang2020a} for unconstrained or locally set constrained games, \cite{Koshal2016,Liang2017,Deng2019} for aggregative games, and \cite{Yi2019,Lu2019,Pavel2020} for generalized games.
Recently, NE seeking agorithms for $N$-cluster games have started to draw researchers' attention \cite{Ye2018,Ye2019,Ye2020,Zeng2019,Nian2021,Sun2021a}. Specifically, the work in \cite{Ye2018} proposed a NE seeking algorithm based on a dynamic average consensus and the gradient play, which was extended in \cite{Ye2019} to reduce the communication and computation costs by introducing an interference graph among the agents in the same cluster. Then, the $N$-cluster game was solved by an extremum seeking-based approach in \cite{Ye2020}. Different from the full-decision information setup in \cite{Ye2018,Ye2019,Ye2020}, the work in \cite{Zeng2019} modeled the decisions of the agents in the same cluster by a decision vector that needs to be agreed on, and introduced a leader-following hierarchy inter-cluster communication mechanism, where only the cluster leader can exchange information across different clusters. A distributed NE seeking algorithm based on subgradient dynamics was designed. The works in \cite{Nian2021,Sun2021a} also considered the partial-decision information setup, and supposed that all agents from all clusters are connected. Distributed NE seeking algorithms based on gradient-play and primal-dual dynamics were presented, respectively. The work in \cite{Tatarenko2021} considered the same cluster-level decision vector modeling as in \cite{Zeng2019}, and developed a discrete-time NE seeking algorithm based on gradient play with no explicit communications between clusters. As compared to the primal methods based on gradient play, gradient tracking is preferred for its faster convergence speed without sacrificing the accuracy, and hence has found its great advantages in the fields of $N$-cluster games recently \cite{Pang2022,Pang2022a,Meng2020,Zimmermann2021}. Specifically, the works in \cite{Pang2022,Pang2022a} studied an $N$-cluster game with full-decision information setup, and proposed gradient-free NE seeking algorithms based on gradient tracking. In the consideration of partial-decision information game setup, the work in \cite{Meng2020} considered the same cluster-level decision vector modeling and applied gradient tracking with the leader-following hierarchy inter-cluster communication mechanism as in \cite{Zeng2019}. The work in \cite{Zimmermann2021} also considered such cluster-level decision vector modeling, and proposed a gradient tracking based method with a leaderless inter-cluster communication mechanism, where each agent maintains an estimate of all cluster-level decisions. In fact, non-cooperative games have the nature that each agent (player) has only authority to its own decision variable and may only be able to control its own action in the practical applications. Though being theoretically equivalent, it would be natural to model the decision variable of each agent individually \cite{Ye2018,Ye2019,Ye2020,Sun2021a,Pang2022,Pang2022a}, as compared to modeling the collection of all agents as a single decision vector \cite{Zeng2019,Nian2021,Meng2020,Zimmermann2021,Tatarenko2021}. For discrete-time methods, step-size is usually assumed to be publicly known by all agents, and hence a common step-size agreed by all agents is chosen in the algorithm development \cite{Tatarenko2021,Pang2022,Meng2020,Zimmermann2021}. Concerning a distributed setup, it would be preferred that agents are allowed to select their own preferred step-sizes. However, introducing such uncoordinated step-sizes may not always guarantee the convergence, and also brings substantial complexity in the convergence analysis. For example, the work in \cite{Pang2022a} allows agents from different clusters to use different step-sizes and proved the convergence to the NE under certain conditions. However, agents from the same cluster still adopt a common step-size.

\textit{Contributions:} 
In view of the comparison with the existing literature, the main contributions of this paper are summarized as follows.
1). As compared to modeling the collection of all agents as a single decision vector \cite{Zeng2019,Nian2021,Meng2020,Zimmermann2021,Tatarenko2021}, this paper models the decision variable of each agent in each cluster individually, which is more ameanable to the applications where agents can only control their own actions. Moreover, we consider a partial-decision information game setup, \textit{i.e.}, the agents only have direct access to their own decisions. Such game setup makes it applicable to the scenarios where there is no central node to collect and broadcast the agents' information, and has a wider range of applications. 
2). We propose a distributed NE seeking algorithm for the $N$-cluster game based on gradient tracking. To account for the partial-decision information game setup, the agents estimate other agents' decisions via a leader-following consensus protocol across clusters over a high-level network. This is different from the leader-following hierarchy inter-cluster communication mechanism as in \cite{Zeng2019,Meng2020} where there is a pre-defined leader to exchange information across clusters. Meanwhile, the agents in the same cluster perform the gradient tracking updates over a low-level network. Moreover, we consider directed graphs for both high-level and low-level networks, and their corresponding adjacency matrices are not required to be doubly stochastic. Hence, the matrix configurations are simplified, since not all directed graphs admit corresponding doubly stochastic adjacency matrices \cite{Gharesifard2010}. 
3). Different from all existing discrete-time methods for $N$-cluster games where either a common step-size is publicly known by all agents \cite{Tatarenko2021,Pang2022,Meng2020,Zimmermann2021} or only known by agents from the same cluster \cite{Pang2020a}, we consider fully uncoordinated step-sizes, \textit{i.e.}, all agents (both among and across) the clusters are allowed to choose their own preferred constant step-sizes, which makes the algorithm even more distributed.

\textit{Notations:}
We use $\mathbf{1}_m$ for an $m$-dimensional vector with all elements being 1, and $I_m$ for an $m\times m$ identity matrix. For a vector $\pi$, we use $\text{diag}(\pi)$ to denote a diagonal matrix formed by the elements of $\pi$. For any two vectors $u, v$, their inner product is denoted by $\langle u, v \rangle$; 
their weighted inner product due to a positive vector (\textit{i.e.}, vector with all elements being positive) $\pi$ is denoted by $\langle u,v\rangle_\pi \triangleq u^\top \text{diag}(\pi)^{-1} v$.
The transpose of $u$ is denoted by $u^\top$. Moreover, we use $\|u\|_2$ for its standard Euclidean norm, \textit{i.e.}, $\|u\|_2 = \sqrt{\langle u, u\rangle}$, and $\|u\|_\pi$ for its weighted Euclidean norm due to $\pi$, \textit{i.e.}, $\|u\|_\pi = \|\text{diag}(\sqrt{\pi})^{-1}u\|_2$.
For vector $a$, we use $[a]_i$ to denote its $i$-th entry. 
The transpose and spectral norm of a matrix $A$ are denoted by $A^\top$ and $\|A\|_2$, respectively. The matrix norm $\|A\|_\pi$ induced by $\|\cdot\|_\pi$ is defined as $\|A\|_\pi \triangleq \|\text{diag}(\sqrt{\pi})^{-1}A\text{ diag}(\sqrt{\pi})\|_2$. 
We use $\rho(A)$ to represent the spectral radius of a square matrix $A$, and $A_\infty$ to indicate its infinite power (if it exists) $\lim_{k\to\infty}A^k$.

\section{Problem Statement}\label{sec:problem_formulation}

An $N$-cluster game, defined by $\Gamma(\mathcal{N},\{J^i\},\{\mathbb{R}^{n_i}\})$, is a multi-player non-cooperative game played among $N$ clusters, where each cluster, indexed by $i\in\mathcal{N} \triangleq \{1, 2, \ldots, N\}$, consists of a group of agents, denoted by $\mathcal{V}^i \triangleq \{1,2,\ldots,n_i\}$, to cooperatively minimize a cluster-level cost function $J^i$. Denote $n \triangleq \sum_{i=1}^Nn_i$. Then, the cluster-level cost function $J^i:\mathbb{R}^n\to\mathbb{R}$ is defined as 
\begin{equation*}
J^i(\mathbf{x}) \triangleq \frac1{n_i}\sum_{j=1}^{n_i}J^i_j(\mathbf{x}^i,\mathbf{x}^{-i}), \quad \forall i\in\mathcal{N},
\end{equation*}
where $J^i_j(\mathbf{x})$ is a local cost function of agent $j$ in cluster $i$, $\mathbf{x}^i \triangleq [x^{i\top}_1,\ldots,x^{i\top}_{n_i}]^\top \in \mathbb{R}^{n_i}$ is a collection of all agents' decisions in cluster $i$ with $x^i_j\in\mathbb{R}$ being the action of agent $j$ in cluster $i$, $\mathbf{x}^{-i}\in \mathbb{R}^{n-n_i}$ denotes a collection of all agents' decisions except cluster $i$, and $\mathbf{x}\triangleq [\mathbf{x}^{1\top},\ldots,\mathbf{x}^{N\top}]^\top$.

\begin{Definition}
(NE of $N$-Cluster Games). A vector $\mathbf{x}^* \triangleq (\mathbf{x}^{i*},\mathbf{x}^{-i*})\in\mathbb{R}^n$ is said to be an NE of the $N$-cluster non-cooperative game $\Gamma(\mathcal{N},\{J^i\},\{\mathbb{R}^{n_i}\})$, if and only if
$J^i(\mathbf{x}^{i*},\mathbf{x}^{-i*})\leq J^i(\mathbf{x}^i,\mathbf{x}^{-i*})$, $\forall i\in \mathcal{N}$.
\end{Definition} 

Motivated by the work in \cite{Ye2017c}, we assume the agents are equipped with two level networks: a high-level network for agents' decisions exchange across clusters and $N$ low-level networks for agents' (partial) gradient exchange within the cluster. 
For the low-level network in each cluster $i\in\mathcal{N}$, it is a directed graph consisting of all agents in the same cluster, denoted by $\mathcal{G}_i(\mathcal{V}^i,\mathcal{E}^i)$ with an adjacency matrix $A^i\triangleq[a^i_{jk}]\in\mathbb{R}^{n_i\times n_i}$, $a^i_{jk}>0$ if $(k,j)\in\mathcal{E}^i$ and $a^i_{jk}=0$ otherwise. We assume $(k,k)\in\mathcal{E}^i,\forall k\in\mathcal{V}^i$. 
For the high-level network, it is a directed graph consisting of all agents in all clusters, denoted by $\bar{\mathcal{G}}(\bar{\mathcal{V}},\bar{\mathcal{E}})$ with an adjacency matrix $\bar{A}\triangleq [\bar{a}_{pq}] \in \mathbb{R}^{n\times n}$, $\bar{a}_{pq}>0$ if $(q,p)\in\bar{\mathcal{E}}$ and $\bar{a}_{pq}=0$ otherwise. We assume $(p,p)\in\bar{\mathcal{E}},\forall p\in\bar{\mathcal{V}}$. 
The following two standard assumptions on the two level networks are imposed.

\begin{Assumption}\label{assumption_graph}
For $i\in \mathcal{N}$, the digraph $\mathcal{G}_i$ is strongly connected. The associated adjacency matrix $A^i$ is column stochastic, \textit{i.e.}, $\mathbf{1}_{n_i}^\top A^i = \mathbf{1}_{n_i}^\top$.
\end{Assumption}
\begin{Assumption}\label{assumption_graph_high_level}
The digraph $\bar{\mathcal{G}}$ is strongly connected. Its associated adjacency matrix $\bar{A}$ is row stochastic, \textit{i.e.}, $\bar{A}\mathbf{1}_{n} = \mathbf{1}_{n}$.
\end{Assumption}

Under Assumption~\ref{assumption_graph}, it is known that $A^i$ is primitive and column stochastic, then we denote its right eigenvector corresponding to the eigenvalue of 1 by $\pi^i\triangleq [\pi^i_1,\ldots,\pi^i_{n_i}]^\top$, such that $\mathbf{1}_{n_i}^\top\pi^i = 1$. Then, $\pi^i$ corresponds to $A^i$'s non-$\mathbf{1}_{n_i}$ Perron vector with eigenvalue 1, and hence all elements in $\pi^i$ are positive, and $A^i_\infty = \pi^i\mathbf{1}_{n_i}^\top$. Define $\iota \triangleq \max_{i\in\mathcal{N}}\|\pi^i\|_2$, $\bm{\pi}^i \triangleq n_i\pi^i$, and $\bm{\pi} \triangleq [\bm{\pi}^{1\top},\ldots,\bm{\pi}^{N\top}]^\top$. Denote the smallest and largest elements of $\bm{\pi}$ by $\underline{\bm{\pi}}$ and $\overline{\bm{\pi}}$, respectively. With the above notations, we can obtain that $\sqrt{\underline{\bm{\pi}}}\|\cdot\|_{\bm{\pi}^i}\leq\|\cdot\|_2\leq \sqrt{\overline{\bm{\pi}}}\|\cdot\|_{\bm{\pi}^i}$ and $\sqrt{\underline{\bm{\pi}}}\|\cdot\|_{\bm{\pi}}\leq\|\cdot\|_2\leq \sqrt{\overline{\bm{\pi}}}\|\cdot\|_{\bm{\pi}}$ based on the definitions of the weighted Euclidean norm, which will be frequently applied in the subsequent analysis.

Then, the following contraction properties of the adjacency matrix $A^i$ can be deduced.
\begin{Lemma}\label{lemma:property_A}
(see \cite[Lemma~1]{Xin2019}) Under Assumption~\ref{assumption_graph}, the adjacency matrix $A^i$ holds that
$\|A^i - A^i_\infty\|_{\bm{\pi}^i} < 1$, and $\|I_{n_i} - A^i_\infty\|_{\bm{\pi}^i} =1$.
\end{Lemma}

Next, we make the following assumption on the agents' local cost functions.

\begin{Assumption}\label{assumption_local_f_lipschitz}
For each $j\in\mathcal{V}^i,i\in\mathcal{N}$, the local cost function $J^i_j(\mathbf{x}^i,\mathbf{x}^{-i})$ is convex, continuously differentiable in $\mathbf{x}^i$; the partial gradient with respect to $x^i_k$, $\forall k \in \mathcal{V}^i$ (denoted by $\nabla^i_kJ^i_j(\mathbf{x})$ for simplicity), is $\mathcal{L}$-Lipschitz continuous in $\mathbf{x}$, \textit{i.e.}, for any $\mathbf{x},\mathbf{x}'\in\mathbb{R}^n$, we have $\|\nabla^i_kJ^i_j(\mathbf{x}) - \nabla^i_kJ^i_j(\mathbf{x}')\|_2\leq \mathcal{L}\|\mathbf{x}- \mathbf{x}'\|_2$. 
\end{Assumption}

The game mapping of $\Gamma(\mathcal{N},\{J^i\},\{\mathbb{R}^{n_i}\})$ is defined as
$\bm{F}(\mathbf{x})\triangleq [\nabla_{\mathbf{x}^1}J^1(\mathbf{x})^\top,\ldots,\nabla_{\mathbf{x}^N}J^N(\mathbf{x})^\top]^\top$.
Then, it follows from Assumption~\ref{assumption_local_f_lipschitz} that $\bm{F}(\mathbf{x})$ is Lipschitz continuous, \textit{i.e.}, for any $\mathbf{x},\mathbf{x}'\in\mathbb{R}^n$, we have $\|\bm{F}(\mathbf{x}) - \bm{F}(\mathbf{x}')\|_2 \leq \sqrt{n}\mathcal{L}\|\mathbf{x}- \mathbf{x}'\|_2$. Next, the following assumption on the game mapping condition is supposed.

\begin{Assumption}\label{assumption_game_mapping}
The game mapping $\bm{F}$ of game $\Gamma$ is strongly monotone on $\mathbb{R}^n$ with a constant $\chi>0$, \textit{i.e.}, for any $\mathbf{x},\mathbf{x}'\in\mathbb{R}^n$, we have $\langle \bm{F}(\mathbf{x})-\bm{F}(\mathbf{x}'), \mathbf{x}-\mathbf{x}' \rangle\geq \chi\|\mathbf{x}-\mathbf{x}'\|_2^2$.
\end{Assumption}
\begin{Remark}
Under Assumptions~\ref{assumption_local_f_lipschitz} and \ref{assumption_game_mapping}, game $\Gamma$ admits a unique NE $\mathbf{x}^*$. Moreover, at NE, $\bm{F}(\mathbf{x}^*) = \mathbf{0}_n$, and hence $\langle \bm{F}(\mathbf{x}), \mathbf{x}-\mathbf{x}^* \rangle\geq \chi\|\mathbf{x}-\mathbf{x}^*\|_2^2$.
\end{Remark}

\section{Algorithm}\label{sec:distr_opt}

In this section, we present a distributed NE seeking strategy for the $N$-cluster game under partial-decision information scenario, followed by the detailed convergence analysis. 

For the notational convenience, agent $j\in\mathcal{V}^i$ in cluster $i\in\mathcal{N}$ of the low-level network is referred to agent $\sum_{l=0}^{i-1}n_l+j$ with $n_0 = 0$ in the high-level network. Hence, its action variable $x^i_j$ is relabeled by $\mathsf{y}_{\sum_{l=0}^{i-1}n_l+j}$, \textit{i.e.}, $x^i_j=\mathsf{y}_{\sum_{l=0}^{i-1}n_l+j}$. Then, each agent $j\in\mathcal{V}^i$, $i\in\mathcal{N}$ needs to maintain the action variable $x^i_j$, and gradient tracker variables $g^i_{jk}$ for $\forall k\in\mathcal{V}^i$. Moreover, in the high-level network, each agent $p\in\bar{\mathcal{V}}$ also needs to maintain an estimation variable $y^p_q$ for the action $\mathsf{y}_q$ of agent $q\in\bar{\mathcal{V}}$. For $i\in\mathcal{N}$, $p\in\bar{\mathcal{V}}$, we denote that 
$\mathbf{y}^p_i \triangleq [y^p_{\sum_{l=0}^{i-1}n_l+1},\ldots,y^p_{\sum_{l=0}^{i-1}n_l+n_i}]^\top\in\mathbb{R}^{n_i}, 
\mathbf{y}^p \triangleq [\mathbf{y}^{p\top}_1,\ldots,\mathbf{y}^{p\top}_N]^\top\in\mathbb{R}^{n},
\mathbf{y}^{(i)} \triangleq [\mathbf{y}^{\sum_{l=0}^{i-1}n_l+1\top},\ldots,\mathbf{y}^{\sum_{l=0}^{i-1}n_l+n_i\top}]^\top\in\mathbb{R}^{n_in}$.
We use subscript $t$ to denote the values of all these variables at time-step $t$. The update laws are designed as follows.
\begin{subequations}\label{eq:algorithm_distr}
\begin{align}
\shortintertext{high-level network:}
&y^p_{q,t+1} = \sum_{l=1}^n\bar{a}_{pl}y^l_{q,t}+\delta_p\bar{a}_{pq}(\mathsf{y}_{q,t} - y^p_{q,t}), \quad p,q\in\bar{\mathcal{V}}\label{eq:update_est_distr}
\shortintertext{low-level network:}
&x^i_{j,t+1} = x^i_{j,t}-\gamma^i_j g^i_{jj,t},\quad j\in\mathcal{V}^i,i\in\mathcal{N}\label{eq:update_x_distr}\\
&g^i_{jk,t+1} = \sum_{l=1}^{n_i}a^i_{jl} g^i_{lk,t} + \nabla^i_kJ^i_j(\mathbf{y}^{\sum_{l=0}^{i-1}n_l+j}_{t+1})\nonumber\\
&\quad\quad\quad\quad - \nabla^i_kJ^i_j(\mathbf{y}^{\sum_{l=0}^{i-1}n_l+j}_t),\quad j,k\in\mathcal{V}^i,i\in\mathcal{N}\label{eq:update_y_distr}
\end{align}
\end{subequations}
with arbitrary $x^i_{j,0}\in\mathbb{R}$, $\mathbf{y}^{\sum_{l=0}^{i-1}n_l+j}_0 \in \mathbb{R}^n$ and $g^i_{jk,0} = \nabla^i_kJ^i_j(\mathbf{y}^{\sum_{l=0}^{i-1}n_l+j}_0)$,
where 
$\delta_p$ is a constant parameter for agent $p\in\bar{\mathcal{V}}$ to be determined later, and $\gamma^i_j > 0$ is a constant step-size sequence adopted by agent $j\in\mathcal{V}^i$, $i\in\mathcal{N}$. 
Denote the largest step-size by ${\gamma}_{M}\triangleq\max_{j\in\mathcal{V}^i,i\in\mathcal{N}}\gamma^i_j$ and the average of all step-sizes by $\bar{\gamma} \triangleq \frac1n\sum_{j\in\mathcal{V}^i,i\in\mathcal{N}}\gamma^i_j$.
Define the heterogeneity of the step-size as the following ratio, $\epsilon_\gamma\triangleq \|\bm{\gamma}-\bar{\bm{\gamma}}\|_2/\|\bar{\bm{\gamma}}\|_2$, where $\bm{\gamma} \triangleq [\gamma^1_1,\ldots,\gamma^1_{n_1},\ldots,\gamma^N_1,\ldots,\gamma^N_{n_N}]^\top$ and $\bar{\bm{\gamma}} \triangleq \bar{\gamma}\mathbf{1}_n$.  


\section{Convergence Analysis} \label{sec:conv_analysis}

The following notations are made throughout the convergence analysis for convenience. For $\forall k\in\mathcal{V}^i, i \in\mathcal{N}$,
$\mathbf{g}^i_{k,t} \triangleq [g^i_{1k,t},\ldots,g^i_{n_ik,t}]^\top\in\mathbb{R}^{n_i}$, 
$\bar{\mathbf{g}}^i_{k,t} \triangleq \frac1{n_i}\mathbf{1}_{n_i}^\top\mathbf{g}^i_{k,t}\in\mathbb{R}$, 
$\mathbf{g}_t \triangleq [g^1_{11,t},g^1_{22,t},\ldots,g^N_{n_Nn_N,t}]^\top\in\mathbb{R}^n$,
$\nabla^i_k \mathbf{J}^i(\mathbf{y}^{(i)}_t) \triangleq [\nabla^i_k J^i_1(\mathbf{y}^{\sum_{l=0}^{i-1}n_l+1}_t),\ldots,\nabla^i_k J^i_{n_i}(\mathbf{y}^{\sum_{l=0}^{i-1}n_l+n_i}_t)]^\top\in\mathbb{R}^{n_i}$, and
$\nabla^i_k \bar{\mathbf{J}}^i(\mathbf{y}^{(i)}_t) \triangleq \frac1{n_i}\mathbf{1}_{n_i}^\top\nabla^i_k \mathbf{J}^i(\mathbf{y}^{(i)}_t)\in\mathbb{R}$.
Then, the concatenated form of \eqref{eq:algorithm_distr} is given by 
\begin{subequations}
\begin{align}
\mathbf{y}_{q,t+1} &= \bar{A}\mathbf{y}_{q,t} \nonumber\\
&\quad+ \text{diag}([\delta_1\bar{a}_{1q},\ldots,\delta_n\bar{a}_{nq}]^\top)(\mathbf{1}_n\mathsf{y}_{q,t} - \mathbf{y}_{q,t})\label{eq:update_est_distr_compact}\\
\mathbf{x}_{t+1} &= \mathbf{x}_t-\text{diag}(\bm{\gamma}) \mathbf{g}_t,\label{eq:update_x_distr_compact}\\
\mathbf{g}^i_{k,t+1} &= A^i \mathbf{g}^i_{k,t} + \nabla^i_k \mathbf{J}^i(\mathbf{y}^{(i)}_{t+1}) - \nabla^i_k \mathbf{J}^i(\mathbf{y}^{(i)}_t).\label{eq:update_y_distr_compact}
\end{align}
\end{subequations}

The convergence analysis of the proposed algorithm is conducted by establishing a linear system, which is composed of three major expressions: i) the total decision estimation error $\sum_{q=1}^n\|\mathbf{y}_{q,t}-\mathbf{1}_n\mathsf{y}_{q,t}\|^2_e$ where $\|\cdot\|_e$ norm is introduced later, ii) the total gradient tracking error $\sum_{i=1}^N\sum_{k=1}^{n_i}\|\mathbf{g}^i_{k,t} - A^i_\infty\mathbf{g}^i_{k,t}\|_{\bm{\pi}^i}^2$, and iii) the gap between all agents' decisions and the NE $\|\mathbf{x}_t-\mathbf{x}^*\|_{\bm{\pi}}^2$.

\subsection{Auxiliary Results}

We first start with an important property of the adjacency matrix $\bar{A}$ in the following lemma.
\begin{Lemma}\label{lemma:A_matrix_high_level}
Under Assumption~\ref{assumption_graph_high_level}, let $\delta_p>0, p \in \bar{\mathcal{V}}$ be chosen such that $0\leq \delta_p\bar{a}_{pq} < 2\bar{a}_{pp}$ $\forall q\in\bar{\mathcal{V}}$. Then, the matrix $\tilde{A}_q\triangleq [\tilde{a}^q_{pm}], q\in\bar{\mathcal{V}}$ with its entry given by
\begin{align*}
\tilde{a}^q_{pm} = \begin{cases}\bar{a}_{pm} & \text{if } p\neq m \\|\bar{a}_{pp} -\delta_p\bar{a}_{pq}|& \text{if } p = m\end{cases}
\end{align*}
holds that $\rho(\tilde{A}_q)<1$. Moreover, there exists a matrix norm $\|\cdot\|_E$ such that $\|\tilde{A}_q\|_E < 1$ for $\forall q\in\bar{\mathcal{V}}$.
\end{Lemma}
\begin{Proof}
See Appendix~A for the proof.
\end{Proof}

We denote the vector norm which is compatible with the matrix norm $\|\cdot\|_E$ by $\|\cdot\|_e$, \textit{i.e.}, $\|Av\|_e \leq \|A\|_E\|v\|_e$ for a matrix $A$ and a vector $v$ with compatible size. Due to the equivalence of all norms in a finite-dimensional vector space, there exists $\underline{C}>0$ and $\overline{C}>0$ such that $\underline{C}\|v\|_2 \leq \|v\|_e \leq \overline{C}\|v\|_2$.

Next, we provide a bound on the stacked gradient estimator $\mathbf{g}_t$ in the following lemma.

\begin{Lemma}\label{lemma:stacked_grad_tracker_distr}
Under Assumptions~\ref{assumption_graph} and \ref{assumption_local_f_lipschitz}, the stacked gradient tracker $\mathbf{g}_t$ holds that
\begin{align*}
\|\mathbf{g}_t\|_2^2&\leq3\overline{\bm{\pi}}\iota^2\mathcal{L}^2\sum_{i=1}^Nn_i^3\|\mathbf{x}_t - \mathbf{x}^*\|_{\bm{\pi}}^2+ 3\overline{\bm{\pi}}\sum_{i=1}^N\sum_{k=1}^{n_i}\|\mathbf{g}^i_{k,t} \\
&\quad- A^i_\infty\mathbf{g}^i_{k,t}\|_{\bm{\pi}^i}^2 + \frac{3n^2\iota^2\mathcal{L}^2}{\underline{C}^2}\sum_{q=1}^n\|\mathbf{y}_{q,t}-\mathbf{1}_n\mathsf{y}_{q,t}\|^2_e.
\end{align*}
\end{Lemma}
\begin{Proof}
For $k\in\mathcal{V}^i, i \in\mathcal{N}$,
\begin{align*}
&\|\mathbf{g}^i_{k,t}\|_2 \leq \|\mathbf{g}^i_{k,t} - \mathbf{v}^i\mathbf{1}_{n_i}^\top\mathbf{g}^i_{k,t}\|_2 + \|\mathbf{v}^i\mathbf{1}_{n_i}^\top\mathbf{g}^i_{k,t}\|_2\\
&\quad\leq \sqrt{\overline{\bm{\pi}}}\|\mathbf{g}^i_{k,t} - A^i_\infty\mathbf{g}^i_{k,t}\|_{\bm{\pi}^i} + n_i\iota\|\bar{\mathbf{g}}^i_{k,t} - \nabla^i_k \bar{\mathbf{J}}^i(\mathbf{y}^{(i)}_t)\|_2 \\
&\quad\quad+ n_i\iota\|\nabla^i_k \bar{\mathbf{J}}^i(\mathbf{y}^{(i)}_t) - \nabla^i_k J^i(\mathbf{x}_t)\|_2 \\
&\quad\quad+ n_i\iota\|\nabla^i_k J^i(\mathbf{x}_t) - \nabla^i_k J^i(\mathbf{x}^*)\|_2,
\end{align*}
where $\nabla^i_k J^i(\mathbf{x}^*) = 0$ has been applied.
Since $A^i$ is column stochastic, it follows from \eqref{eq:update_y_distr_compact} that
$\bar{\mathbf{g}}^i_{k,t+1} = \bar{\mathbf{g}}^i_{k,t} + \nabla^i_k \bar{\mathbf{J}}^i(\mathbf{y}^{(i)}_{t+1}) - \nabla^i_k \bar{\mathbf{J}}^i(\mathbf{y}^{(i)}_t)$.
Due to the initial conditions $\mathbf{g}^i_{k,0} = \nabla^i_k \mathbf{J}^i(\mathbf{y}^{(i)}_0)$, we obtain
\begin{align}\label{eq:averaged_grad_tracker_distr}
\bar{\mathbf{g}}^i_{k,t} = \nabla^i_k \bar{\mathbf{J}}^i(\mathbf{y}^{(i)}_t).
\end{align}
Thus,
$\|\mathbf{g}^i_{k,t}\|_2^2 \leq 3\overline{\bm{\pi}}\|\mathbf{g}^i_{k,t} - A^i_\infty\mathbf{g}^i_{k,t}\|_{\bm{\pi}^i}^2
+ 3n_i^2\iota^2\|\nabla^i_k \bar{\mathbf{J}}^i(\mathbf{y}^{(i)}_t) - \nabla^i_k J^i(\mathbf{x}_t)\|_2^2
+ 3\overline{\bm{\pi}}n_i^2\iota^2\mathcal{L}^2\|\mathbf{x}_t - \mathbf{x}^*\|_{\bm{\pi}}^2$.
Since
$\|\mathbf{g}_t\|_2^2= \sum_{i=1}^N\sum_{k=1}^{n_i}\|g^i_{kk,t}\|_2^2\leq\sum_{i=1}^N\sum_{k=1}^{n_i}\|\mathbf{g}^i_{k,t}\|_2^2$,
substituting the above relation, and noting that
\begin{align}
&\sum_{i=1}^N\sum_{k=1}^{n_i}n_i^2\|\nabla^i_k \bar{\mathbf{J}}^i(\mathbf{y}^{(i)}_t) - \nabla^i_k J^i(\mathbf{x}_t)\|^2_2\nonumber\\
&\quad\leq\sum_{i=1}^N\sum_{k=1}^{n_i}\sum_{j=1}^{n_i}n_i\|\nabla^i_k J^i_j(\mathbf{y}^{\sum_{l=0}^{i-1}n_l+j}_t)-\nabla^i_k J^i_j(\mathbf{x}_t)\|^2_2\nonumber\\
&\quad\leq\mathcal{L}^2\sum_{i=1}^N\sum_{k=1}^{n_i}\sum_{j=1}^{n_i}n_i\|\mathbf{y}^{\sum_{l=0}^{i-1}n_l+j}_t-\mathbf{x}_t\|^2_2\nonumber\\
&\quad=\mathcal{L}^2\sum_{i=1}^N\sum_{k=1}^{n_i}\sum_{j=1}^{n_i}\sum_{q=1}^nn_i\|y^{\sum_{l=0}^{i-1}n_l+j}_{q,t}-\mathsf{y}_{q,t}\|^2_2\nonumber\\
&\quad\leq \frac{n^2\mathcal{L}^2}{\underline{C}^2}\sum_{q=1}^n\|\mathbf{y}_{q,t}-\mathbf{1}_n\mathsf{y}_{q,t}\|^2_e,\label{eq:second_difference_term}
\end{align}
we obtain the desired result.
\end{Proof}

Next, we establish the inequality iterations of the three major terms in Lemmas~\ref{lemma:action_estimation_error}, \ref{lemma:optimality_gap_distr} and \ref{lemma:grad_track_error_distr}, respectively. 

We first derive a bound on the total action estimation error characterized by $\sum_{q=1}^n\|\mathbf{1}_n\mathsf{y}_{q,t} - \mathbf{y}_{q,t}\|_e^2$.
\begin{Lemma}\label{lemma:action_estimation_error}
Under Assumptions~\ref{assumption_graph}, \ref{assumption_graph_high_level} and \ref{assumption_local_f_lipschitz}, the total action estimation error satisfies:
\begin{align*}
&\sum_{q=1}^n\|\mathbf{1}_n\mathsf{y}_{q,t+1} - \mathbf{y}_{q,t+1}\|_e^2\\
&\leq \bigg[\frac{1+\bar{\sigma}^2_2}2+\frac{3n^3\varsigma_2\overline{C}^2\iota^2\mathcal{L}^2{\gamma}_{M}^2}{\underline{C}^2}\bigg]\sum_{q=1}^n\|\mathbf{1}_n\mathsf{y}_{q,t} - \mathbf{y}_{q,t}\|_e^2\\
&\quad + 3\overline{\bm{\pi}}n\varsigma_2\overline{C}^2{\gamma}_{M}^2 \sum_{i=1}^N\sum_{k=1}^{n_i}\|\mathbf{g}^i_{k,t} - A^i_\infty\mathbf{g}^i_{k,t}\|_{\bm{\pi}^i}^2\\
&\quad +3\overline{\bm{\pi}}n\varsigma_2\overline{C}^2\iota^2\mathcal{L}^2\sum_{i=1}^Nn_i^3{\gamma}_{M}^2\|\mathbf{x}_t - \mathbf{x}^*\|_{\bm{\pi}}^2.
\end{align*}
where $\bar{\sigma}_2$, $\varsigma_2$ and $\sigma_{\tilde{A}_q}$ are some constants.
\end{Lemma}
\begin{Proof}
From \eqref{eq:update_x_distr}, we have for $q \in \bar{\mathcal{V}}$,
$\mathbf{1}_n\mathsf{y}_{q,t+1} = \mathbf{1}_n\mathsf{y}_{q,t} - \gamma_q\mathbf{1}_n\mathbf{g}_{q,t}$,
where $\gamma_q$ and $\mathbf{g}_{q,t}$ denote the step-size and gradient tracker of agent $q$, respectively. That is, $\gamma^i_j$ and $g^i_{jj,t}$ correspond to $\gamma_q$ and $\mathbf{g}_{q,t}$ with $q = \sum_{l=0}^{i-1}n_l+j$, if agent $j\in\mathcal{V}^i$ in cluster $i\in\mathcal{N}$ is considered.
Subtracting it by \eqref{eq:update_est_distr_compact} and taking the vector norm $\|\cdot\|_e$ on both sides, we have
\begin{align*}
\|\mathbf{1}_n\mathsf{y}_{q,t+1} &- \mathbf{y}_{q,t+1}\|_e\leq \|(\bar{A}-\text{diag}([\delta_1\bar{a}_{1q},\ldots,\delta_n\bar{a}_{nq}]^\top))\\
&\quad\quad(\mathbf{1}_n\mathsf{y}_{q,t} - \mathbf{y}_{q,t}) - \gamma_q \mathbf{1}_n\mathbf{g}_{q,t}\|_e\\
&\leq\|\tilde{A}_q(\mathbf{1}_n\mathsf{y}_{q,t} - \mathbf{y}_{q,t})\|_e + \gamma_q \|\mathbf{1}_n\mathbf{g}_{q,t}\|_e\\
&\leq\|\tilde{A}_q\|_E\|\mathbf{1}_n\mathsf{y}_{q,t} - \mathbf{y}_{q,t}\|_e + \sqrt{n}\overline{C}{\gamma}_{M} \|\mathbf{g}_{q,t}\|_2.
\end{align*}
Define $\sigma_{\tilde{A}_q}\triangleq \|\tilde{A}_q\|_E$ and square both sides, we obtain
\begin{align}
&\|\mathbf{1}_n\mathsf{y}_{q,t+1} - \mathbf{y}_{q,t+1}\|^2_e\leq \sigma_{\tilde{A}_q}^2\|\mathbf{1}_n\mathsf{y}_{q,t} - \mathbf{y}_{q,t}\|^2_e + n{\gamma}_{M}^2\|\mathbf{g}_{q,t}\|^2_e\nonumber\\
&\quad\quad+\frac{1-\sigma_{\tilde{A}_q}^2}2\|\mathbf{1}_n\mathsf{y}_{q,t} - \mathbf{y}_{q,t}\|_e^2 + \frac{2\sigma_{\tilde{A}_q}^2n\overline{C}^2{\gamma}_{M}^2}{1-\sigma_{\tilde{A}_q}^2}\|\mathbf{g}_{q,t}\|_2^2\nonumber\\
&\quad=\frac{1+\bar{\sigma}_2^2}2\|\mathbf{1}_n\mathsf{y}_{q,t} - \mathbf{y}_{q,t}\|_e^2 + n\varsigma_2\overline{C}^2{\gamma}_{M}^2 \|\mathbf{g}_{q,t}\|_2^2,\label{eq:estimation_error}
\end{align}
where we denote $\bar{\sigma}_2 \triangleq \max_{q\in\bar{\mathcal{V}}}\sigma_{\tilde{A}_q}$, $\varsigma_2\triangleq \max_{q\in\bar{\mathcal{V}}}\frac{1+\sigma_{\tilde{A}_q}^2}{1-\sigma_{\tilde{A}_q}^2}$. Summing over $q = 1$ to $n$ gives
$\sum_{q=1}^n\|\mathbf{1}_n\mathsf{y}_{q,t+1} - \mathbf{y}_{q,t+1}\|_e^2\leq \frac{1+\bar{\sigma}^2_2}2\sum_{q=1}^n\|\mathbf{1}_n\mathsf{y}_{q,t} - \mathbf{y}_{q,t}\|_e^2
+ n\varsigma_2\overline{C}^2{\gamma}_{M}^2 \|\mathbf{g}_t\|_2^2$.
Substituting the result in Lemma~\ref{lemma:stacked_grad_tracker_distr} completes the proof.
\end{Proof}

Next, we bound the gap between all agents' decisions and the NE, characterized by $\|\mathbf{x}_t - \mathbf{x}^*\|_{\bm{\pi}}^2$.
\begin{Lemma}\label{lemma:optimality_gap_distr}
Under Assumptions~\ref{assumption_graph}, \ref{assumption_graph_high_level}, \ref{assumption_local_f_lipschitz} and \ref{assumption_game_mapping}, the agents' decisions $\mathbf{x}_t$ satisfies that
\begin{align*}
&\|\mathbf{x}_{t+1}-\mathbf{x}^*\|_{\bm{\pi}}^2\leq \bigg[\frac{3\overline{\bm{\pi}}{\gamma}_{M}^2}{\underline{\bm{\pi}}}+ \frac{2{\gamma}_{M}^2}{\underline{\bm{\pi}}\chi\bar{\gamma}}\bigg]\sum_{i=1}^N\sum_{k=1}^{n_i}\|\mathbf{g}^i_{k,t} - A^i_\infty\mathbf{g}^i_{k,t}\|_{\bm{\pi}^i}^2\\
&+\bigg[1-\underline{\bm{\pi}}\chi\bar{\gamma}+\frac{3\overline{\bm{\pi}}\iota^2\mathcal{L}^2{\gamma}_{M}^2}{\underline{\bm{\pi}}}\sum_{i=1}^Nn_i^3+2\overline{\bm{\pi}}\sqrt{n}\mathcal{L}\epsilon_\gamma\bar{\gamma}\bigg]\|\mathbf{x}_t - \mathbf{x}^*\|_{\bm{\pi}}^2\\
&\quad+\bigg[\frac{3n^2\iota^2\mathcal{L}^2{\gamma}_{M}^2}{\underline{\bm{\pi}}\underline{C}^2}+ \frac{2\overline{\bm{\pi}}n\mathcal{L}^2{\gamma}_{M}^2}{\underline{\bm{\pi}}\underline{C}^2\chi\bar{\gamma}}\bigg]\sum_{q=1}^n\|\mathbf{y}_{q,t}-\mathbf{1}_n\mathsf{y}_{q,t}\|^2_e.
\end{align*}
\end{Lemma}
\begin{Proof}
It follows from \eqref{eq:update_x_distr_compact} that
$\mathbf{x}_{t+1}-\mathbf{x}^* = \mathbf{x}_t - \text{diag}(\bm{\gamma}) \mathbf{g}_t - \mathbf{x}^*$.
Taking the norm on both sides gives
\begin{subequations}\label{eq:optimality_gap_expand_distr}
\begin{align}
&\|\mathbf{x}_{t+1}-\mathbf{x}^*\|_{\bm{\pi}}^2 = \|\mathbf{x}_t - \text{diag}(\bm{\gamma}) \mathbf{g}_t - \mathbf{x}^*\|_{\bm{\pi}}^2\nonumber\\
&\leq \|\mathbf{x}_t - \mathbf{x}^*\|_{\bm{\pi}}^2 + \frac{{\gamma}_{M}^2}{\underline{\bm{\pi}}}\|\mathbf{g}_t\|_2^2\nonumber\\
&\quad-2\langle \mathbf{x}_t - \mathbf{x}^*, \text{diag}(\bm{\gamma})(\mathbf{g}_t-\text{diag}(\bm{\pi})\bm{F}(\mathbf{x}_t))\rangle_{\bm{\pi}}\label{eq:optimality_gap_expand_1_distr}\\
&\quad-2\langle \mathbf{x}_t - \mathbf{x}^*, \text{diag}(\bm{\gamma}- \bar{\bm{\gamma}})\text{diag}((\bm{\pi})\bm{F}(\mathbf{x}_t)\rangle_{\bm{\pi}}\label{eq:optimality_gap_expand_2_distr}\\
&\quad-2\bar{\gamma}\langle \mathbf{x}_t - \mathbf{x}^*, \text{diag}(\bm{\pi})\bm{F}(\mathbf{x}_t)\rangle_{\bm{\pi}}. \label{eq:optimality_gap_expand_3_distr}
\end{align}
\end{subequations}

For \eqref{eq:optimality_gap_expand_1_distr}, it follows that
\begin{align*}
&-2\langle \mathbf{x}_t - \mathbf{x}^*, \text{diag}(\bm{\gamma})(\mathbf{g}_t - \text{diag}(\bm{\pi})\bm{F}(\mathbf{x}_t))\rangle_{\bm{\pi}} \\
&= -2\sum_{i=1}^N\sum_{k=1}^{n_i}\gamma^i_k\langle x^i_{k,t} - x^{i*}_k,  g^i_{kk,t}- n_i\pi^i_k\nabla^i_k J^i(\mathbf{x}_t)\rangle_{n_i\pi^i_k}\\
&=-2\sum_{i=1}^N\sum_{k=1}^{n_i}\bigg(\gamma^i_k\langle x^i_{k,t} - x^{i*}_k, g^i_{kk,t} - n_i\pi^i_k\bar{\mathbf{g}}^i_{k,t}\rangle_{n_i\pi^i_k}\\
&\quad+\gamma^i_k\langle x^i_{k,t} - x^{i*}_k,  n_i\pi^i_k(\bar{\mathbf{g}}^i_{k,t} - \nabla^i_k \bar{\mathbf{J}}^i(\mathbf{y}^{(i)}_t))\rangle_{n_i\pi^i_k}\\
&\quad+\gamma^i_k\langle x^i_{k,t} - x^{i*}_k, n_i\pi^i_k(\nabla^i_k \bar{\mathbf{J}}^i(\mathbf{y}^{(i)}_t) - \nabla^i_k J^i(\mathbf{x}_t))\rangle_{n_i\pi^i_k}\bigg).
\end{align*}
The first part holds that
\begin{align*}
&-2\sum_{i=1}^N\sum_{k=1}^{n_i}\gamma^i_k\langle x^i_{k,t} - x^{i*}_k, g^i_{kk,t} - n_i\pi^i_k\bar{\mathbf{g}}^i_{k,t}\rangle_{n_i\pi^i_k}\\
&=-2\sum_{i=1}^N\sum_{k=1}^{n_i}\gamma^i_k\langle x^i_{k,t} - x^{i*}_k, g^i_{kk,t} - \pi^i_k\mathbf{1}_{n_i}^\top\mathbf{g}^i_{k,t} \rangle_{n_i\pi^i_k}\\
&\leq2{\gamma}_{M}\sum_{i=1}^N\sum_{k=1}^{n_i}\|x^i_{k,t} - x^{i*}_k\|_{n_i\pi^i_k}\|g^i_{kk,t} - \pi^i_k\mathbf{1}_{n_i}^\top\mathbf{g}^i_{k,t}\|_{n_i\pi^i_k}\\
&\leq\frac{\underline{\bm{\pi}}\chi\bar{\gamma}}2\|\mathbf{x}_t - \mathbf{x}^*\|_{\bm{\pi}}^2+ \frac{2{\gamma}_{M}^2}{\underline{\bm{\pi}}\chi\bar{\gamma}}\sum_{i=1}^N\sum_{k=1}^{n_i}\|\mathbf{g}^i_{k,t} - A^i_\infty\mathbf{g}^i_{k,t}\|_{\bm{\pi}^i}^2.
\end{align*}
For the second part, it follows from \eqref{eq:averaged_grad_tracker_distr} that
\begin{align*}
&\langle x^i_{k,t} - x^{i*}_k, n_i\pi^i_k(\bar{\mathbf{g}}^i_{k,t} - \nabla^i_k \bar{\mathbf{J}}^i(\mathbf{y}^{(i)}_t))\rangle_{n_i\pi^i_k}=0.
\end{align*}
For the third part, we have
\begin{align*}
&-2\sum_{i=1}^N\sum_{k=1}^{n_i}\gamma^i_k\langle x^i_{k,t} - x^{i*}_k, n_i\pi^i_k(\nabla^i_k \bar{\mathbf{J}}^i(\mathbf{y}^{(i)}_t) \\
&\quad\quad - \nabla^i_k J^i(\mathbf{x}_t))\rangle_{n_i\pi^i_k}\leq2{\gamma}_{M}\sum_{i=1}^N\sum_{k=1}^{n_i}\| x^i_{k,t} - x^{i*}_k\|_{n_i\pi^i_k} \\
&\quad\quad\quad\times\sqrt{n_i\pi^i_k}\|\nabla^i_k \bar{\mathbf{J}}^i(\mathbf{y}^{(i)}_t) - \nabla^i_k J^i(\mathbf{x}_t)\|_2\\
&\quad\leq\frac{\underline{\bm{\pi}}\chi\bar{\gamma}}2\|\mathbf{x}_t - \mathbf{x}^*\|_{\bm{\pi}}^2 \\
&\quad\quad+ \frac{2\overline{\bm{\pi}}{\gamma}_{M}^2}{\underline{\bm{\pi}}\chi\bar{\gamma}}\sum_{i=1}^N\sum_{k=1}^{n_i}n_i\|\nabla^i_k \bar{\mathbf{J}}^i(\mathbf{y}^{(i)}_t) - \nabla^i_k J^i(\mathbf{x}_t)\|^2_2.
\end{align*}
The last term follows the same derivation as in \eqref{eq:second_difference_term} that
$\sum_{i=1}^N\sum_{k=1}^{n_i}n_i\|\nabla^i_k \bar{\mathbf{J}}^i(\mathbf{y}^{(i)}_t) - \nabla^i_k J^i(\mathbf{x}_t)\|^2_2
\leq \frac{n\mathcal{L}^2}{\underline{C}^2}\sum_{q=1}^n\|\mathbf{y}_{q,t}-\mathbf{1}_n\mathsf{y}_{q,t}\|^2_e$.
Hence, the third part can be further obtained that
\begin{align*}
&-2\sum_{i=1}^N\sum_{k=1}^{n_i}\gamma^i_k\langle x^i_{k,t} - x^{i*}_k,  \\
&\quad\quad\quad\quad n_i\pi^i_k(\nabla^i_k \bar{\mathbf{J}}^i(\mathbf{y}^{(i)}_t) - \nabla^i_k J^i(\mathbf{x}_t))\rangle_{n_i\pi^i_k}\\
&\quad\leq\frac{\underline{\bm{\pi}}\chi\bar{\gamma}}2\|\mathbf{x}_t - \mathbf{x}^*\|_{\bm{\pi}}^2 + \frac{2\overline{\bm{\pi}}n\mathcal{L}^2{\gamma}_{M}^2}{\underline{\bm{\pi}}\underline{C}^2\chi\bar{\gamma}} \sum_{q=1}^n\|\mathbf{y}_{q,t}-\mathbf{1}_n\mathsf{y}_{q,t}\|^2_e.
\end{align*}
Combining the above three parts, we obtain that
\begin{align}
&-2\bar{\gamma}\langle \mathbf{x}_t - \mathbf{x}^*, \mathbf{g}_t - \text{diag}(\bm{\pi})\bm{F}(\mathbf{x}_t)\rangle_{\bm{\pi}}\nonumber\\
&\quad\leq \underline{\bm{\pi}}\chi\bar{\gamma}\|\mathbf{x}_t - \mathbf{x}^*\|_{\bm{\pi}}^2 + \frac{2{\gamma}_{M}^2}{\underline{\bm{\pi}}\chi\bar{\gamma}}\sum_{i=1}^N\sum_{k=1}^{n_i}\|\mathbf{g}^i_{k,t} - A^i_\infty\mathbf{g}^i_{k,t}\|_{\bm{\pi}^i}^2\nonumber\\
&\quad+ \frac{2\overline{\bm{\pi}}n\mathcal{L}^2{\gamma}_{M}^2}{\underline{\bm{\pi}}\underline{C}^2\chi\bar{\gamma}} \sum_{q=1}^n\|\mathbf{y}_{q,t}-\mathbf{1}_n\mathsf{y}_{q,t}\|^2_e.\label{eq:optimality_gap_expand_1_result_distr}
\end{align}

For \eqref{eq:optimality_gap_expand_2_distr}, it follows that
\begin{align}
&-2\langle \mathbf{x}_t - \mathbf{x}^*, \text{diag}(\bm{\gamma}- \bar{\bm{\gamma}})\text{diag}((\bm{\pi})\bm{F}(\mathbf{x}_t)\rangle_{\bm{\pi}} \nonumber\\
&\quad=-2\langle \mathbf{x}_t - \mathbf{x}^*, \text{diag}(\bm{\gamma}- \bar{\bm{\gamma}})\text{diag}((\bm{\pi})(\bm{F}(\mathbf{x}_t)-\bm{F}(\mathbf{x}^*))\rangle_{\bm{\pi}} \nonumber\\
&\quad\leq 2\|\mathbf{x}_t - \mathbf{x}^*\|_{\bm{\pi}}\|\text{diag}(\bm{\gamma}- \bar{\bm{\gamma}})\text{diag}((\bm{\pi})(\bm{F}(\mathbf{x}_t)-\bm{F}(\mathbf{x}^*))\|_{\bm{\pi}} \nonumber\\
&\quad\leq2\sqrt{\overline{\bm{\pi}}}\|\mathbf{x}_t - \mathbf{x}^*\|_{\bm{\pi}}\|\text{diag}(\bm{\gamma}- \bar{\bm{\gamma}})(\bm{F}(\mathbf{x}_t)-\bm{F}(\mathbf{x}^*))\|_2 \nonumber\\
&\quad\leq 2\overline{\bm{\pi}}\sqrt{n}\mathcal{L}\epsilon_\gamma\bar{\gamma}\|\mathbf{x}_t-\mathbf{x}^*\|_{\bm{\pi}}^2, \label{eq:optimality_gap_expand_2_result_distr}
\end{align}
since $\|\text{diag}(\bm{\gamma} - \bar{\bm{\gamma}})(\bm{F}(\mathbf{x}_t)-\bm{F}(\mathbf{x}^*))\|_2\leq \|\bm{\gamma} - \bar{\bm{\gamma}}\|_2\|\bm{F}(\mathbf{x}_t)-\bm{F}(\mathbf{x}^*)\|_2\leq \sqrt{n}\mathcal{L}\epsilon_\gamma\bar{\gamma}\|\mathbf{x}_t-\mathbf{x}^*\|_2\leq\sqrt{\overline{\bm{\pi}}}\sqrt{n}\mathcal{L}\epsilon_\gamma\bar{\gamma}\|\mathbf{x}_t-\mathbf{x}^*\|_{\bm{\pi}}$.

For \eqref{eq:optimality_gap_expand_3_distr}, by Assumption~\ref{assumption_game_mapping}, we have
\begin{align}
&-2\bar{\gamma}\langle \mathbf{x}_t - \mathbf{x}^*, \text{diag}(\bm{\pi})\bm{F}(\mathbf{x}_t)\rangle_{\bm{\pi}}=-2\bar{\gamma}\langle \mathbf{x}_t - \mathbf{x}^*, \bm{F}(\mathbf{x}_t)\rangle\nonumber\\
&\quad\leq-2\chi\bar{\gamma}\|\mathbf{x}_t - \mathbf{x}^*\|_2^2\leq-2\underline{\bm{\pi}}\chi\bar{\gamma}\|\mathbf{x}_t - \mathbf{x}^*\|_{\bm{\pi}}^2.\label{eq:optimality_gap_expand_3_result_distr}
\end{align}

Substituting \eqref{eq:optimality_gap_expand_1_result_distr}, \eqref{eq:optimality_gap_expand_2_result_distr}, \eqref{eq:optimality_gap_expand_3_result_distr} and Lemma~\ref{lemma:stacked_grad_tracker_distr} into \eqref{eq:optimality_gap_expand_distr} yields the desired result.
\end{Proof}

Finally, we quantify the total gradient tracking error, measured by $\sum_{i=1}^N\sum_{k=1}^{n_i}\|\mathbf{g}^i_{k,t} - A^i_\infty\mathbf{g}^i_{k,t}\|_{\bm{\pi}^i}^2$.

\begin{Lemma}\label{lemma:grad_track_error_distr}
Under Assumptions~\ref{assumption_graph}, \ref{assumption_graph_high_level} and \ref{assumption_local_f_lipschitz}, the total gradient tracking error $\sum_{i=1}^N\sum_{k=1}^{n_i}\|\mathbf{g}^i_{k,t} - A^i_\infty\mathbf{g}^i_{k,t}\|_{\bm{\pi}^i}^2$ satisfies
\begin{align*}
&\sum_{i=1}^N\sum_{k=1}^{n_i}\|\mathbf{g}^i_{k,t+1} - A^i_\infty\mathbf{g}^i_{k,t+1}\|_{\bm{\pi}^i}^2 \leq\bigg[\frac{1+\bar{\sigma}^2_1}2\\
&\quad+\frac{9\overline{\bm{\pi}}n^2\varsigma_1(1+n\varsigma_2\overline{C}^2)\mathcal{L}^2{\gamma}_{M}^2}{\underline{\bm{\pi}}}\bigg]\sum_{i=1}^N\sum_{k=1}^{n_i} \|\mathbf{g}^i_{k,t} - A^i_\infty\mathbf{g}^i_{k,t}\|_{\bm{\pi}^i}^2\\
&+\bigg[\frac{9n^4\varsigma_1(1+n\varsigma_2\overline{C}^2)\iota^2\mathcal{L}^4{\gamma}_{M}^2}{\underline{\bm{\pi}}\underline{C}^2}\\
&\quad\quad\quad\quad\quad+\frac{3n\varsigma_1(3+\bar{\sigma}_2^2)\mathcal{L}^2}{2\underline{\bm{\pi}}\underline{C}^2}\bigg]\sum_{q=1}^n\|\mathbf{y}_{q,t}-\mathbf{1}_n\mathsf{y}_{q,t}\|_e^2\\
&+\frac{9\overline{\bm{\pi}}n^2\varsigma_1(1+n\varsigma_2\overline{C}^2)\iota^2\mathcal{L}^4{\gamma}_{M}^2}{\underline{\bm{\pi}}}\sum_{i=1}^Nn_i^3\|\mathbf{x}_t - \mathbf{x}^*\|_{\bm{\pi}}^2,
\end{align*}
where $\bar{\sigma}_1$, $\varsigma_1$ and $\sigma_{A^i}$ are some constants.
\end{Lemma}
\begin{Proof}
It is obtained from \eqref{eq:update_y_distr_compact} that
\begin{align*}
&\|\mathbf{g}^i_{k,t+1} - A^i_\infty\mathbf{g}^i_{k,t+1}\|_{\bm{\pi}^i}^2= \|A^i\mathbf{g}^i_{k,t} - A^i_\infty\mathbf{g}^i_{k,t}\|_{\bm{\pi}^i}^2 + \|(I_{n_i}\\
&-A^i_\infty)(\nabla^i_k \mathbf{J}^i(\mathbf{y}^{(i)}_{t+1}) - \nabla^i_k \mathbf{J}^i(\mathbf{y}^{(i)}_t))\|_{\bm{\pi}^i}^2 +2\langle A^i\mathbf{g}^i_{k,t} \\
&- A^i_\infty\mathbf{g}^i_{k,t}, (I_{n_i}-A^i_\infty)(\nabla^i_k \mathbf{J}^i(\mathbf{y}^{(i)}_{t+1}) - \nabla^i_k \mathbf{J}^i(\mathbf{y}^{(i)}_t) )\rangle_{\bm{\pi}^i}.
\end{align*}
Define $\sigma_{A^i} \triangleq \|A^i-A^i_\infty\|_{\bm{\pi}^i}$. 
It is noted that $\|I_{n_i}-A^i_\infty\|_{\bm{\pi}^i} = 1$ from Lemma~\ref{lemma:property_A}, then
\begin{align}
&\|\mathbf{g}^i_{k,t+1} - A^i_\infty\mathbf{g}^i_{k,t+1}\|_{\bm{\pi}^i}^2 \nonumber\\
&\quad\leq \sigma_{A_i}^2\|\mathbf{g}^i_{k,t} - A^i_\infty\mathbf{g}^i_{k,t}\|_{\bm{\pi}^i}^2 + \|\nabla^i_k \mathbf{J}^i(\mathbf{y}^{(i)}_{t+1})\nonumber\\
&\quad\quad - \nabla^i_k \mathbf{J}^i(\mathbf{y}^{(i)}_t)\|_{\bm{\pi}^i}^2+2\| A^i\mathbf{g}^i_{k,t} \nonumber\\
&\quad\quad- A^i_\infty\mathbf{g}^i_{k,t}\|_{\bm{\pi}^i} \|\nabla^i_k \mathbf{J}^i(\mathbf{y}^{(i)}_{t+1}) - \nabla^i_k \mathbf{J}^i(\mathbf{y}^{(i)}_t)\|_{\bm{\pi}^i}\nonumber\\
&\quad\leq \sigma_{A_i}^2\|\mathbf{g}^i_{k,t} - A^i_\infty\mathbf{g}^i_{k,t}\|_{\bm{\pi}^i}^2 \nonumber\\
&\quad\quad+ \|\nabla^i_k \mathbf{J}^i(\mathbf{y}^{(i)}_{t+1}) - \nabla^i_k \mathbf{J}^i(\mathbf{y}^{(i)}_t)\|_{\bm{\pi}^i}^2\nonumber\\
&\quad\quad +\frac{1-\sigma_{A_i}^2}2\|\mathbf{g}^i_{k,t} - A^i_\infty\mathbf{g}^i_{k,t}\|_{\bm{\pi}^i}^2\nonumber\\
&\quad\quad+\frac{2\sigma_{A_i}^2}{1-\sigma_{A_i}^2}\|\nabla^i_k \mathbf{J}^i(\mathbf{y}^{(i)}_{t+1}) - \nabla^i_k \mathbf{J}^i(\mathbf{y}^{(i)}_t)\|_{\bm{\pi}^i}^2\nonumber\\
&\quad\leq\frac{1+\bar{\sigma}_1^2}2\|\mathbf{g}^i_{k,t} - A^i_\infty\mathbf{g}^i_{k,t}\|_{\bm{\pi}^i}^2\nonumber\\
&\quad\quad+\varsigma_1\|\nabla^i_k \mathbf{J}^i(\mathbf{y}^{(i)}_{t+1}) - \nabla^i_k \mathbf{J}^i(\mathbf{y}^{(i)}_t)\|_{\bm{\pi}^i}^2,\label{eq:grad_track_error_term123_distr}
\end{align}
where $\bar{\sigma}_1 \triangleq \max_{i\in\mathcal{N}}\sigma_{A^i}$ and $\varsigma_1 \triangleq \max_{i\in\mathcal{N}}\frac{(1+\sigma_{A^i}^2)}{1-\sigma_{A^i}^2}$. It follows from Assumption~\ref{assumption_local_f_lipschitz} that
$\|\nabla^i_k \mathbf{J}^i(\mathbf{y}^{(i)}_{t+1}) - \nabla^i_k \mathbf{J}^i(\mathbf{y}^{(i)}_t)\|_{\bm{\pi}^i}^2
\leq \frac{\mathcal{L}^2}{\underline{\bm{\pi}}}\sum_{j=1}^{n_i}\|\mathbf{y}^{\sum_{l=0}^{i-1}n_l+j}_{t+1}-\mathbf{y}^{\sum_{l=0}^{i-1}n_l+j}_t\|_2^2
=\frac{\mathcal{L}^2}{\underline{\bm{\pi}}}\sum_{q=1}^n\sum_{j=1}^{n_i}\|y^{\sum_{l=0}^{i-1}n_l+j}_{q,t+1}-y^{\sum_{l=0}^{i-1}n_l+j}_{q,t}\|_2^2$.
Hence,
\begin{align*}
&\sum_{i=1}^N\sum_{k=1}^{n_i}\|\nabla^i_k \mathbf{J}^i(\mathbf{y}^{(i)}_{t+1}) - \nabla^i_k \mathbf{J}^i(\mathbf{y}^{(i)}_t)\|_{\bm{\pi}^i}^2\\
&\quad\leq\frac{n\mathcal{L}^2}{\underline{\bm{\pi}}}\sum_{q=1}^n\|\mathbf{y}_{q,t+1}-\mathbf{y}_{q,t}\|_2^2\leq\frac{3n\mathcal{L}^2}{\underline{\bm{\pi}}}\sum_{q=1}^n\\
&\quad\quad\bigg(\|\mathbf{y}_{q,t+1}-\mathbf{1}_n\mathsf{y}_{q,t+1}\|_2^2 +\|\mathbf{y}_{q,t}-\mathbf{1}_n\mathsf{y}_{q,t}\|_2^2\\
&\quad\quad+\|\mathbf{1}_n\mathsf{y}_{q,t+1}-\mathbf{1}_n\mathsf{y}_{q,t}\|_2^2\bigg)\\
&\quad\leq\frac{3n\mathcal{L}^2}{\underline{\bm{\pi}}}\sum_{q=1}^n\bigg(\frac{3+\bar{\sigma}_2^2}{2\underline{C}^2}\|\mathbf{y}_{q,t}-\mathbf{1}_n\mathsf{y}_{q,t}\|_e^2 \\
&\quad\quad+\|\mathbf{1}_n\mathsf{y}_{q,t+1}-\mathbf{1}_n\mathsf{y}_{q,t}\|_2^2+ n^2\varsigma_2\overline{C}^2{\gamma}_{M}^2 \|\mathbf{g}_{q,t}\|^2_2\bigg),
\end{align*}
where the last inequality follows from \eqref{eq:estimation_error}. It is noted that
$\sum_{q=1}^n\|\mathbf{1}_n\mathsf{y}_{q,t+1}-\mathbf{1}_n\mathsf{y}_{q,t}\|_2^2 = n\|\mathbf{x}_{t+1}-\mathbf{x}_t\|_2^2
=n\|\mathbf{x}_t -\text{diag}(\bm{\gamma})\mathbf{g}_t-\mathbf{x}_t\|_2^2\leq n{\gamma}_{M}^2\|\mathbf{g}_t\|_2^2$,
Thus, we have
\begin{align*}
&\sum_{i=1}^N\sum_{k=1}^{n_i}\|\nabla^i_k \mathbf{J}^i(\mathbf{y}^{(i)}_{t+1}) - \nabla^i_k \mathbf{J}^i(\mathbf{y}^{(i)}_t)\|_{\bm{\pi}^i}^2\\
&\quad\leq\frac{3n(3+\bar{\sigma}_2^2)\mathcal{L}^2}{2\underline{\bm{\pi}}\underline{C}^2}\sum_{q=1}^n\|\mathbf{y}_{q,t}-\mathbf{1}_n\mathsf{y}_{q,t}\|_e^2\\
&\quad\quad+\frac{3n^2(1+n\varsigma_2\overline{C}^2)\mathcal{L}^2{\gamma}_{M}^2}{\underline{\bm{\pi}}}\|\mathbf{g}_t\|_2^2.
\end{align*}
Summing over $k = 1$ to $n_i$, $i = 1$ to $N$ for \eqref{eq:grad_track_error_term123_distr} and substituting the above result gives
\begin{align*}
&\sum_{i=1}^N\sum_{k=1}^{n_i}\|\mathbf{g}^i_{k,t+1} - A^i_\infty\mathbf{g}^i_{k,t+1}\|_{\bm{\pi}^i}^2 \leq\frac{1+\bar{\sigma}^2_1}2\sum_{i=1}^N\sum_{k=1}^{n_i}\\
&\quad \|\mathbf{g}^i_{k,t}- A^i_\infty\mathbf{g}^i_{k,t}\|_{\bm{\pi}^i}^2 +\frac{3n\varsigma_1(3+\bar{\sigma}_2^2)\mathcal{L}^2}{2\underline{\bm{\pi}}\underline{C}^2}\sum_{q=1}^n\|\mathbf{y}_{q,t}-\mathbf{1}_n\mathsf{y}_{q,t}\|_e^2\\
&\quad+\frac{3n^2\varsigma_1(1+n\varsigma_2\overline{C}^2)\mathcal{L}^2{\gamma}_{M}^2}{\underline{\bm{\pi}}}\|\mathbf{g}_t\|_2^2,
\end{align*}
The proof is completed by substituting Lemma~\ref{lemma:stacked_grad_tracker_distr}.
\end{Proof}

\subsection{Main Results}

Now, we are ready for the analysis on the convergence of the proposed algorithm. With the inequality iterations derived in Lemmas~\ref{lemma:action_estimation_error}, \ref{lemma:optimality_gap_distr} and \ref{lemma:grad_track_error_distr}, we can establish the following linear dynamical system
\begin{align}
\mathbf{u}_{t+1} \leq \mathbf{T}\mathbf{u}_t, \label{eq:dynamical_system_distr}
\end{align}
where
\begin{align*}
\mathbf{u}_t &\triangleq \begin{bmatrix} \|\mathbf{x}_t-\mathbf{x}^*\|_{\bm{\pi}}^2\\\sum_{i=1}^N\sum_{k=1}^{n_i}\|\mathbf{g}^i_{k,t} - A^i_\infty\mathbf{g}^i_{k,t}\|_{\bm{\pi}^i}^2\\\sum_{q=1}^n\|\mathbf{y}_{q,t}-\mathbf{1}_n\mathsf{y}_{q,t}\|^2_e \end{bmatrix},\\
\mathbf{T} &\triangleq\begin{bmatrix}\begin{aligned} &\begin{smallmatrix}1-k_1\bar{\gamma}+k_2{\gamma}_{M}^2\end{smallmatrix}\\ &\quad\quad\begin{smallmatrix}+k_3\epsilon_{\gamma}\bar{\gamma} \end{smallmatrix}\end{aligned} & \begin{smallmatrix} k_4{\gamma}_{M}^2+k_5{\gamma}_{M}^2/\bar{\gamma} \end{smallmatrix} & \begin{smallmatrix} k_6{\gamma}_{M}^2+k_7{\gamma}_{M}^2/\bar{\gamma}\end{smallmatrix} \\ \begin{smallmatrix}k_8{\gamma}_{M}^2\end{smallmatrix} & \begin{smallmatrix}1-k_{10}+k_9{\gamma}_{M}^2\end{smallmatrix} & \begin{smallmatrix}k_{11}+k_{12}{\gamma}_{M}^2\end{smallmatrix} \\ \begin{smallmatrix}k_{13}{\gamma}_{M}^2\end{smallmatrix} & \begin{smallmatrix}k_{14}{\gamma}_{M}^2\end{smallmatrix} & \begin{smallmatrix}1 - k_{16}+k_{15}{\gamma}_{M}^2\end{smallmatrix}\end{bmatrix},
\end{align*}
$k_1\triangleq \underline{\bm{\pi}}\chi$, $k_2\triangleq \frac{3\overline{\bm{\pi}}\iota^2\mathcal{L}^2}{\underline{\bm{\pi}}}\sum_{i=1}^Nn_i^3$, $k_3\triangleq 2\overline{\bm{\pi}}\sqrt{n}\mathcal{L}$, $k_4\triangleq \frac{3\overline{\bm{\pi}}}{\underline{\bm{\pi}}}$, $k_5\triangleq\frac2{\underline{\bm{\pi}}\chi}$, $k_6\triangleq\frac{3n^2\iota^2\mathcal{L}^2}{\underline{\bm{\pi}}\underline{C}^2}$, $k_7\triangleq\frac{2\overline{\bm{\pi}}n\mathcal{L}^2}{\underline{\bm{\pi}}\underline{C}^2\chi}$, $k_8\triangleq\frac{9\overline{\bm{\pi}}n^2\varsigma_1(1+n\varsigma_2\overline{C}^2)\iota^2\mathcal{L}^4}{\underline{\bm{\pi}}}\sum_{i=1}^Nn_i^3$, $k_9 \triangleq \frac{9\overline{\bm{\pi}}n^2\varsigma_1(1+n\varsigma_2\overline{C}^2)\mathcal{L}^2}{\underline{\bm{\pi}}}\sum_{i=1}^Nn_i^3$, $k_{10} \triangleq \frac{1-\bar{\sigma}_1^2}{2}$, $k_{11} \triangleq \frac{3n\varsigma_1(3+\bar{\sigma}_2^2)\mathcal{L}^2}{2\underline{\bm{\pi}}\underline{C}^2}$, $k_{12}\triangleq \frac{9n^4\varsigma_1(1+n\varsigma_2\overline{C}^2)\iota^2\mathcal{L}^4}{\underline{\bm{\pi}}\underline{C}^2}$, $k_{13}\triangleq 3\overline{\bm{\pi}}n\varsigma_2\overline{C}^2\iota^2\mathcal{L}^2\sum_{i=1}^Nn_i^3$, $k_{14} \triangleq 3\overline{\bm{\pi}}n\varsigma_2\overline{C}^2$, $k_{15}\triangleq \frac{3n^3\varsigma_2\overline{C}^2\iota^2\mathcal{L}^2}{\underline{C}^2}$, $k_{16} \triangleq \frac{1-\bar{\sigma}_2^2}{2}$.

Then, the convergence results of all agents' decisions to the NE $\mathbf{x}^*$ can be established based on the convergence of the dynamical system \eqref{eq:dynamical_system_distr}, which are summarized in the following theorem.

\begin{Theorem}\label{theorem:optimality_distr}
Suppose Assumptions~\ref{assumption_graph}, \ref{assumption_graph_high_level}, \ref{assumption_local_f_lipschitz} and \ref{assumption_game_mapping} hold. Generate the agent's action $\{x^i_{j,t}\}_{t\geq0}$, gradient tracker $\{g^i_{jk,t}\}_{t\geq0}$ and estimation variable $\{y^p_{q,t}\}_{t\geq0}$ by \eqref{eq:algorithm_distr} with the uncoordinated constant step-size $\gamma^i_j$ satisfying
\begin{align*}
0\leq\epsilon_\gamma< \frac{k_1}{k_3}, 0<{\gamma}_{M}< \min\bigg\{\frac1{k_1},\gamma^*_1, \gamma^*_2, \gamma^*_3 \bigg\},
\end{align*}
where $\gamma^*_1, \gamma^*_2$ and $\gamma^*_3$ are some constants related to the heterogeneity $\epsilon_\gamma$.
Then, we have $\rho(\mathbf{T})<1$, and $\sup_{\ell\geq t}\|\mathbf{x}_\ell-\mathbf{x}^*\|_{\bm{\pi}}^2$ (respectively, $\sup_{\ell\geq t}\sum_{i=1}^N\sum_{k=1}^{n_i}\|\mathbf{g}^i_{k,\ell} - A^i_\infty\mathbf{g}^i_{k,\ell}\|_{\bm{\pi}^i}^2$ and $\sup_{\ell\geq t}\sum_{q=1}^n\|\mathbf{1}_n\mathsf{y}_{q,\ell} - \mathbf{y}_{q,\ell}\|^2_e$) linearly converges to 0 at a rate of $\rho(\mathbf{T})$.
\end{Theorem}
\begin{Proof}
For system \eqref{eq:dynamical_system_distr}, if $\rho(\mathbf{T})<1$, then $\mathbf{T}^t$ converges to $\mathbf{0}$ at a geometric rate with exponent $\rho(\mathbf{T})$ \cite{Horn1990}, which implies that $\sup_{\ell\geq t}\|\mathbf{x}_\ell-\mathbf{x}^*\|_{\bm{\pi}}^2$, $\sup_{\ell\geq t}\sum_{i=1}^N\sum_{k=1}^{n_i}\|\mathbf{g}^i_{k,\ell} - A^i_\infty\mathbf{g}^i_{k,\ell}\|_{\bm{\pi}^i}^2$ and $\sup_{\ell\geq t}\sum_{q=1}^n\|\mathbf{1}_n\mathsf{y}_{q,\ell} - \mathbf{y}_{q,\ell}\|^2_e$, respectively, converge to 0 with the same rate. 
The following lemma provides a sufficient condition to quantify the spectral radius of a non-negative matrix.
\begin{Lemma}\label{lemma:matrix_spectral_radius}
(see \cite[Cor.~8.1.29]{Horn1990}) Let $A\in\mathbb{R}^{m\times m}$ be a matrix with non-negative entries and $\bm{\theta}\in\mathbb{R}^m$ be a vector with positive entries. If there exists a constant $\lambda\geq0$ such that $A\bm{\theta} < \lambda \bm{\theta}$, then $\rho(A) < \lambda$.
\end{Lemma}

To invoke Lemma~\ref{lemma:matrix_spectral_radius}, the matrix $\mathbf{T}$ has to be non-negative. Thus, it suffices to have
\begin{align}\label{eq:ensure_nonnegative}
0< \bar{\gamma} \leq \frac1{k_1}.
\end{align}
According to Lemma~\ref{lemma:matrix_spectral_radius}, to ensure $\rho(\mathbf{T})<1$, one needs to seek for some positive vector $\bm{\theta} \triangleq [\theta_1,\theta_2,\theta_3]^\top$, where $\theta_1>0$, $\theta_2>0$ and $\theta_3>0$, such that $\mathbf{T}\bm{\theta} < \bm{\theta}$, \textit{i.e.},
\begin{equation*}
\begin{cases}
(1-k_1\bar{\gamma}+k_2{\gamma}_{M}^2+k_3\epsilon_{\gamma}\bar{\gamma})\theta_1 \\
\quad+ (k_4{\gamma}_{M}^2+k_5{\gamma}_{M}^2/\bar{\gamma})\theta_2 +(k_6{\gamma}_{M}^2+k_7{\gamma}_{M}^2/\bar{\gamma})\theta_3< \theta_1,\\
k_8{\gamma}_{M}^2\theta_1+(1-k_{10}+k_9{\gamma}_{M}^2)\theta_2 + (k_{11}+k_{12}{\gamma}_{M}^2)\theta_3 < \theta_2,\\
k_{13}{\gamma}_{M}^2\theta_1+k_{14}{\gamma}_{M}^2\theta_2 + (1-k_{16}+k_{15}{\gamma}_{M}^2)\theta_3 < \theta_3.
\end{cases}
\end{equation*}
Without the loss of generality, we can set $\theta_3 = 1$. Then these inequalities can be further simplified as
\begin{equation}\label{eq:sys_ineq}
\begin{cases}
(k_2\theta_1+k_4\theta_2+k_6)\bar{\gamma}\\
\quad\quad\quad< (k_1-k_3\epsilon_{\gamma})\theta_1\bar{\gamma}^2/{\gamma}_{M}^2 - (k_5\theta_2+k_7),\\
(k_8\theta_1+k_9\theta_2+k_{12}){\gamma}_{M}^2< k_{10}\theta_2-k_{11},\\
(k_{13}\theta_1+k_{14}\theta_2+k_{15}){\gamma}_{M}^2< k_{16}.
\end{cases}
\end{equation}
Therefore, we would like to find the range of the step-size such that \eqref{eq:sys_ineq} hold simultanenously for some $\theta_1>0$, $\theta_2>0$. To ensure the existence of solution ${\gamma}_{M}$ and $\bar{\gamma}$, the right-hand-side of \eqref{eq:sys_ineq} needs to be positive, \textit{i.e.},
$\epsilon_\gamma< \frac{k_1}{k_3}, \theta_1 > \frac{(k_5\theta_2+k_7){\gamma}_{M}^2/\bar{\gamma}^2}{k_1-k_3\epsilon_{\gamma}}, \theta_2 > \frac{k_{11}}{k_{10}}$,
Thus, we can set
$\theta_1 = \frac{2(k_5k_{11}+k_7k_{10}){\gamma}_{M}^2/\bar{\gamma}^2}{k_{10}(k_1-k_3\epsilon_{\gamma})}, \theta_2 = \frac{2k_{11}}{k_{10}}$.
Then, we can solve the three inequalities in \eqref{eq:sys_ineq} respectively. 
For the first inequality, we have
$\bar{\gamma}<[k_7k_{10}(k_1-k_3\epsilon_{\gamma})]/[2k_2(k_5k_{11}+k_7k_{10}){\gamma}_{M}^2/\bar{\gamma}^2+(2k_4k_{11}+k_6k_{10})(k_1-k_3\epsilon_{\gamma})]$.
Noting that $\bar{\gamma}\leq{\gamma}_{M}$ and ${\gamma}_{M}/\bar{\gamma}<n$, for the above inequality to hold, it suffices to have
${\gamma}_{M}< \gamma^*_1$,
where
\begin{align*}
\gamma^*_1 &\triangleq \frac{k_7k_{10}(k_1-k_3\epsilon_\gamma)}{q^*_1},\\
q^*_1 &\triangleq 2n^2k_2(k_5k_{11}+k_7k_{10})+(2k_4k_{11}+k_6k_{10})(k_1-k_3\epsilon_{\gamma}).
\end{align*}

Similarly, for the second and third inequalities to hold, it suffices to have
${\gamma}_{M}< \gamma^*_2,{\gamma}_{M}< \gamma^*_3$,
where
\begin{align*}
\gamma^*_2 &\triangleq \sqrt{\frac{k_{10}k_{11}(k_1-k_3\epsilon_\gamma)}{q^*_2}},
\gamma^*_3 \triangleq \sqrt{\frac{k_{10}k_{16}(k_1-k_3\epsilon_\gamma)}{q^*_3}},\\
q^*_2 &\triangleq 2n^2k_8(k_5k_{11}+k_7k_{10})+(2k_9k_{11}+k_{10}k_{12})(k_1-k_3\epsilon_{\gamma}),\\
q^*_3 &\triangleq 2n^2k_{13}(k_5k_{11}+k_7k_{10})+(2k_{11}k_{14}+k_{10}k_{15})(k_1-k_3\epsilon_{\gamma}).
\end{align*}
Thus, to ensure both \eqref{eq:ensure_nonnegative} and \eqref{eq:sys_ineq} hold simultanenously, it suffices to have
$0\leq\epsilon_\gamma< \frac{k_1}{k_3}, 0<{\gamma}_{M}< \min\{\frac1{k_1},\gamma^*_1, \gamma^*_2, \gamma^*_3 \}$,
which completes the proof.
\end{Proof}
\begin{Remark}
Theorem~\ref{theorem:optimality_distr} shows that the linear convergence of all agents' decisions to the NE is guaranteed when both the largest step-size and its heterogeneity are sufficiently small. Besides, it should be remarked that the bounds on both the largest step-size and heterogeneity are only sufficient but not necessary conditions for the convergence results. The bound conditions could be conservative, and the dynamical system may still converge even though the bound conditions are not satisified.
\end{Remark}

Next, we analyze the convergence of the algorithm when all agents adopt uniform constant step-size, summarized in the following corollary.
\begin{Corollary}\label{corollary:optimality_distr}
Suppose Assumptions~\ref{assumption_graph}, \ref{assumption_graph_high_level}, \ref{assumption_local_f_lipschitz} and \ref{assumption_game_mapping} hold. Generate the agent's action $\{x^i_{j,t}\}_{t\geq0}$, gradient tracker $\{g^i_{jk,t}\}_{t\geq0}$ and estimation variable $\{y^p_{q,t}\}_{t\geq0}$ by \eqref{eq:algorithm_distr} with uniform constant step-size $\gamma$ satisfying
\begin{align*}
0<\gamma< \min\bigg\{\frac1{k_1},\gamma^*_{1,c}, \gamma^*_{2,c}, \gamma^*_{3,c} \bigg\},
\end{align*}
where $\gamma^*_{1,c}, \gamma^*_{2,c}$ and $\gamma^*_{3,c}$ are some constants.
Then, $\sup_{\ell\geq t}\|\mathbf{x}_\ell-\mathbf{x}^*\|_{\bm{\pi}}^2$ (respectively, $\sup_{\ell\geq t}\sum_{i=1}^N\sum_{k=1}^{n_i}\|\mathbf{g}^i_{k,\ell} - A^i_\infty\mathbf{g}^i_{k,\ell}\|_{\bm{\pi}^i}^2$ and $\sup_{\ell\geq t}\sum_{q=1}^n\|\mathbf{1}_n\mathsf{y}_{q,\ell} - \mathbf{y}_{q,\ell}\|^2_e$) linearly converges to 0.
\end{Corollary}
\begin{Proof}
The result directly follows from Theorem~\ref{theorem:optimality_distr} by noting that $\epsilon_\gamma = 0$ and ${\gamma}_{M}/\bar{\gamma}=1$, which gives
\begin{align*}
\gamma^*_{1,c} &\triangleq \frac{k_1k_7k_{10}}{2k_2k_5k_{11}+2k_2k_7k_{10}+2k_1k_4k_{11}+k_1k_6k_{10}},\\
\gamma^*_{2,c} &\triangleq \sqrt{\frac{k_1k_{10}k_{11}}{2k_5k_8k_{11}+2k_7k_8k_{10}+2k_1k_9k_{11}+k_1k_{10}k_{12}}},\\
\gamma^*_{3,c} &\triangleq \sqrt{\frac{k_1k_{10}k_{16}}{2k_5k_{11}k_{13}+2k_7k_{10}k_{13}+2k_1k_{11}k_{14}+k_1k_{10}k_{15}}}.
\end{align*}
Then, following the same arguments in Theorem~\ref{theorem:optimality_distr}, when $0<\gamma<\min\{\frac1{k_1},\gamma^*_{1,c}, \gamma^*_{2,c}, \gamma^*_{3,c}\}$, the convergence of the algorithm is guaranteed.
\end{Proof}

\section{Numerical Simulations}\label{sec:simulation}
\begin{figure}[!tb]
\centering
\subfloat[Low-level 1]{\includegraphics[width=1in]{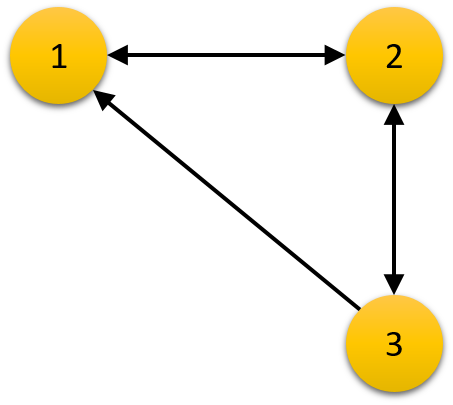}}
\hfil
\subfloat[Low-level 2]{\includegraphics[width=1in]{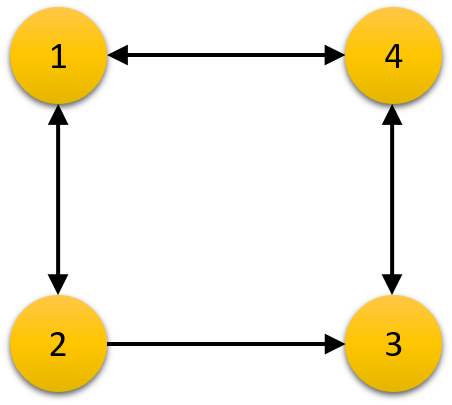}}
\hfil
\subfloat[Low-level 3]{\includegraphics[width=1in]{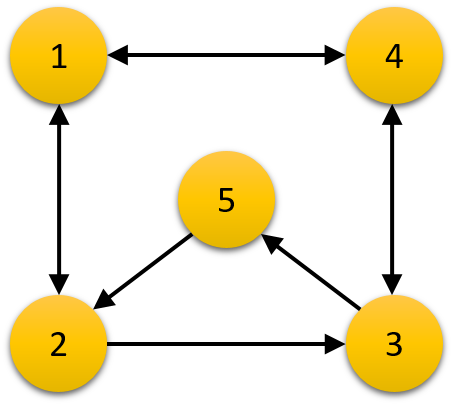}}
\\
\subfloat[High-level network]{\includegraphics[width=3.2in]{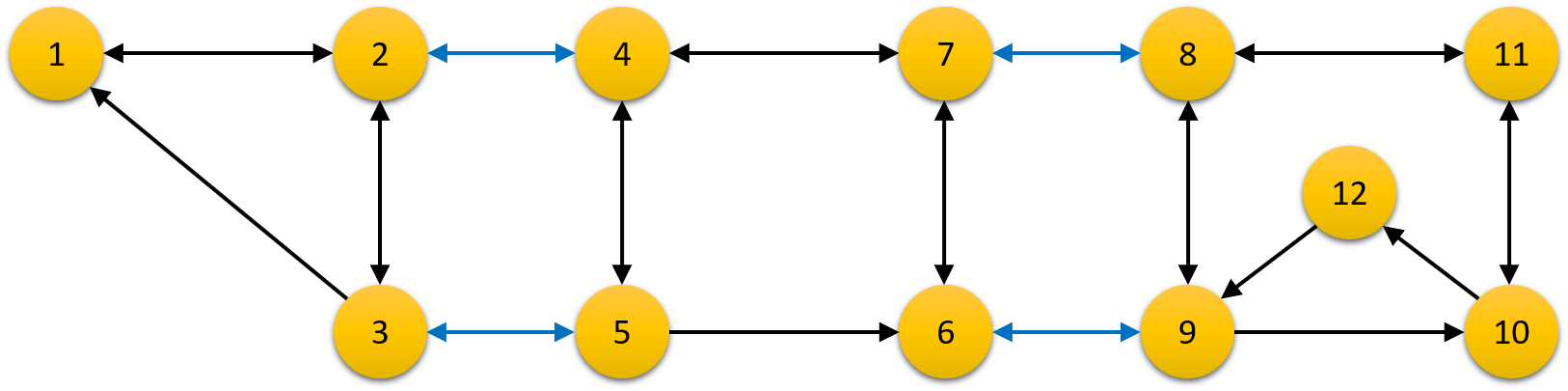}}
\caption{Graph topology of two level networks.}
\label{fig:two_level_networks}
\end{figure}

In this section, we validate the performance of the proposed algorithm by a Cournot competition game. In particular, we consider $N$ firms, and each firm $i\in\mathcal{N}$ consists of $n_i$ branches to help produce goods. For $j\in\mathcal{V}^i, i\in\mathcal{N}$, let $x^i_j$ be the quantity of goods produced by branch $j$ of firm $i$, then its local cost function $J^i_j(\mathbf{x})$ is modeled by the following function
$J^i_j(\mathbf{x}) = c^i_j(x^i_j) - p^i_j(\mathbf{x})x^i_j$,
where $c^i_j(x^i_j) = a^i_j(x^i_j)^2+b^i_j(x^i_j)$ models the cost incurred by generating $x^i_j$ quantity of goods, $p^i_j(\mathbf{x}) = d^i_j - w^{i\top}_j\mathbf{x}$ models the selling price of such goods, $a^i_j, b^i_j, d^i_j\in\mathbb{R}$ and $w^i_j\in\mathbb{R}^n$ are constant parameters.
As a numerical setting, we set $N = 3$, $n_i = 3, 4$ and $5$, respectively. For constant parameters, we let $a^i_j = 1$, $d^i_j = 10+i+j$, $b^i_j$ and each element of $w^i_j$ be uniformly drawn from $[0,1]$, respectively. The two level networks are given in Fig.~\ref{fig:two_level_networks}, which are strongly connected. 
The inital conditions of $\mathbf{x}$ and $\mathbf{y}_q$ are set to some arbitrary values, and $\delta_q = 0.5$, $\forall q\in\bar{\mathcal{V}}$.

\subsection{Algorithm Convergence}\label{subsec:algo_convergence}
\begin{figure}[!tb]
\centering
\subfloat[Cluster 1]{\includegraphics[width=1.7in]{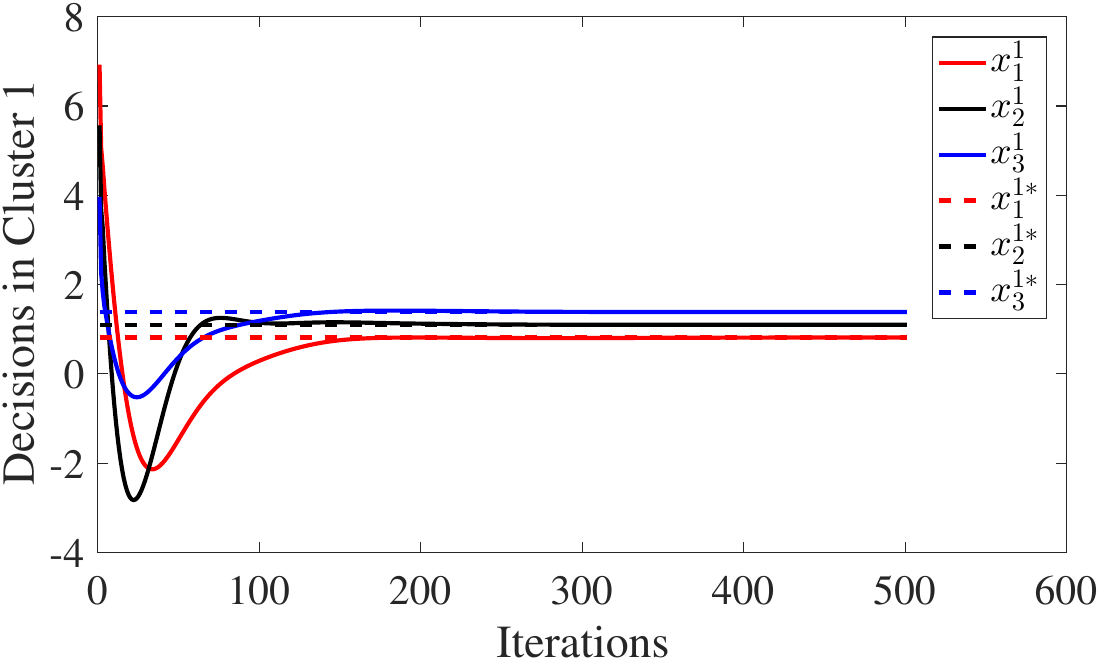}}
\hfil
\subfloat[Cluster 2]{\includegraphics[width=1.7in]{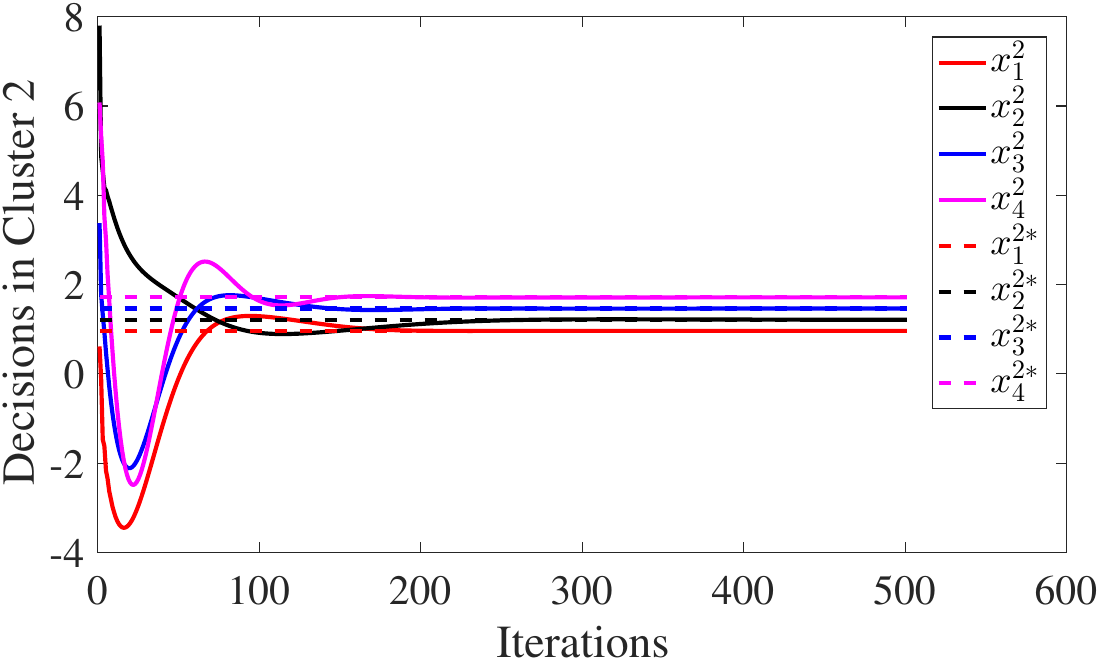}}
\\
\subfloat[Cluster 3]{\includegraphics[width=1.7in]{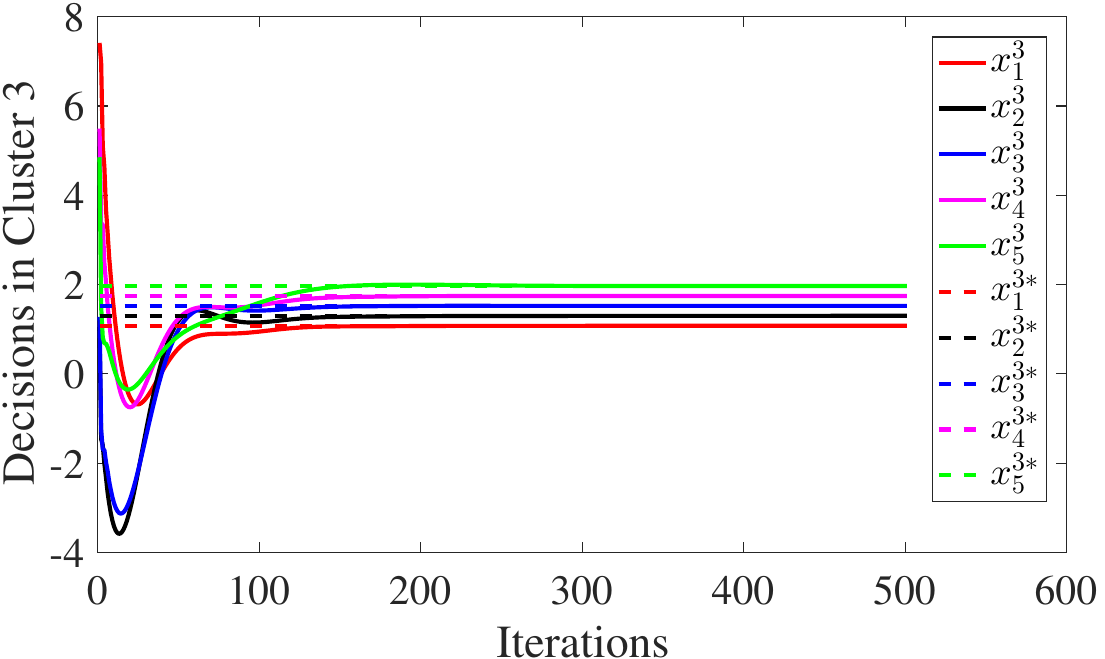}}
\hfil
\subfloat[NE gap]{\includegraphics[width=1.7in]{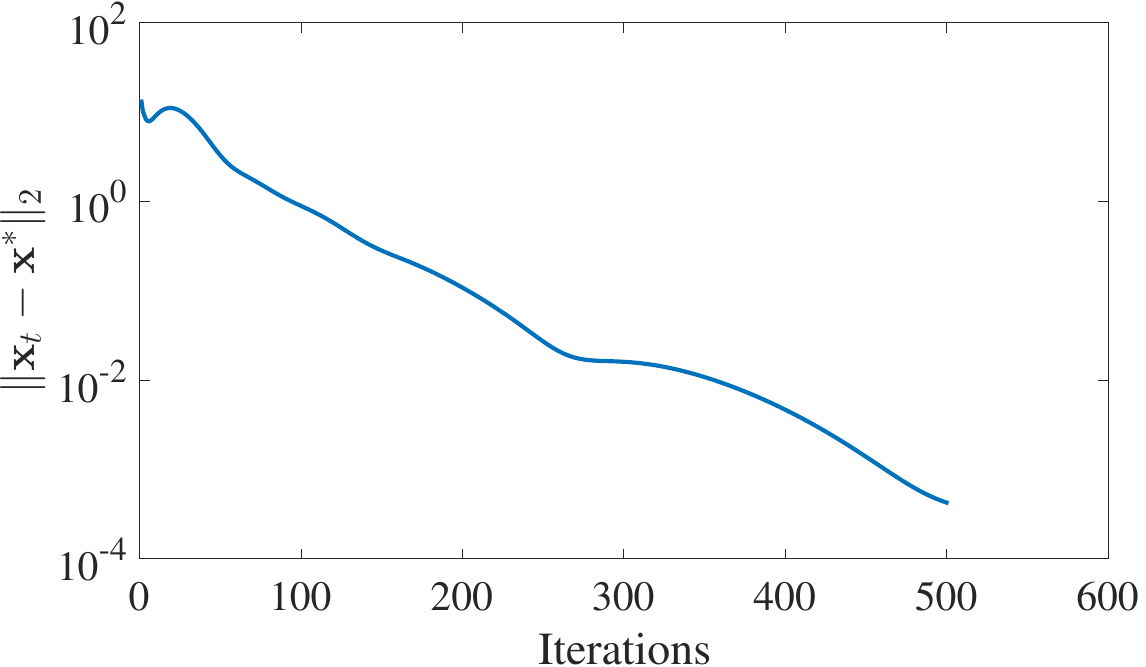}}
\caption{Trajectories of agents' decisions in different clusters.}
\label{fig:action_rgf_con}
\end{figure}

In this part, we focus on the verification of the convergence result derived in Theorem~\ref{theorem:optimality_distr}. The step-size $\gamma^i_j$ is evenly selected from $[0.045,0.1]$, giving a heterogeneity of $0.2381$. Then, the trajectories of the decisions of all firms (and branches) and the NE gap $\|\mathbf{x}_t - \mathbf{x}^*\|_2$ are plotted in Fig.~\ref{fig:action_rgf_con}. As can be seen, the convergence to the NE is obtained and the rate is linear.

\subsection{Influence of step-size on the convergence}
\begin{figure}[!tb]
\centering
\includegraphics[width=2.8in]{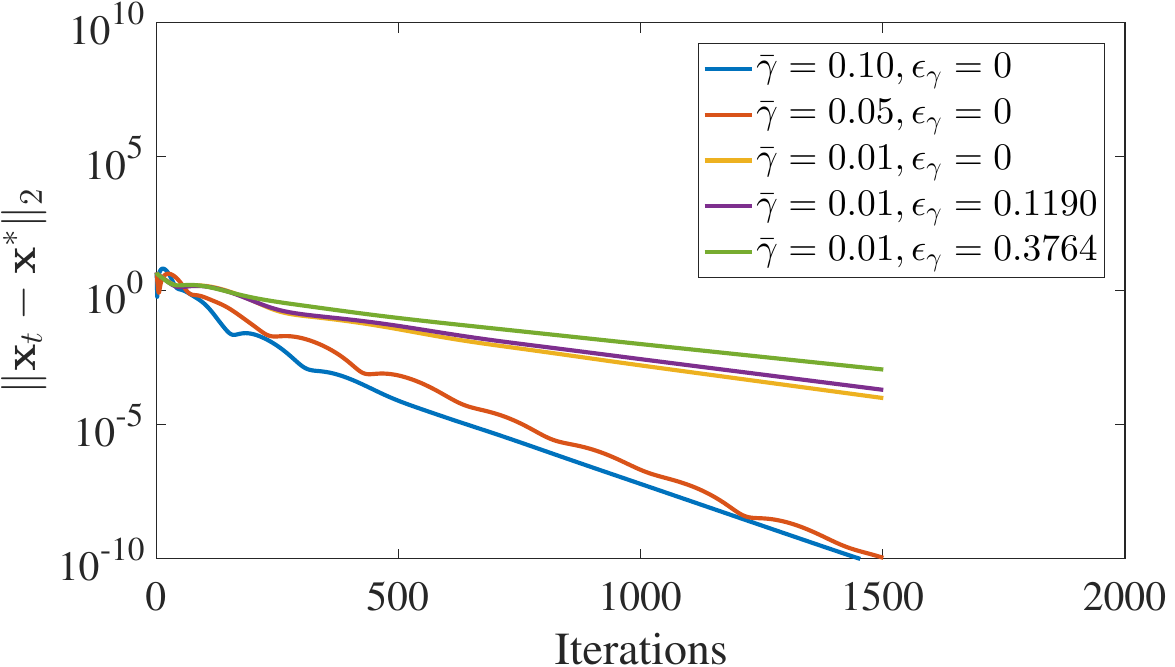} %
\caption{Influence of the step-size and heterogeneity on the rate of convergence.}
\label{fig:influence_step_size}%
\end{figure}

In this part, we investigate the influence of the step-size including the heterogeneity on the convergence. Specifically, we let the agents' step-sizes be selected within $(0,0.1]$. The initial conditions of $\mathbf{x}$ and $\mathbf{y}_q$ are set to zero, while the rest of the parameters are kept the same as in Sec.~\ref{subsec:algo_convergence}. Fig.~\ref{fig:influence_step_size} plots the NE gaps under various step-size cases with different averaged step-size and different heterogeneity. As can be seen, smaller heterogeneity of the step-size and larger averaged step-size lead to a faster rate of convergence.

\section{Conclusions}\label{sec:conclusion}
This paper has considered the $N$-cluster game under partial-decision information settings. A distributed Nash equilibrium (NE) seeking algorithm has been proposed based on a synthesis of leader-following consensus and gradient tracking. It has been shown that all agents' decisions linearly converge to their corresponding NE with uncoordinated constant step-sizes when the largest step-size and the heterogeneity of the step-size are small. The derived result has been validated through a numerical example in a Cournot competition game.



\section*{APPENDIX}

\subsection{Proof of Lemma~\ref{lemma:A_matrix_high_level}}
The first part of the result (\textit{i.e.}, $\rho(\tilde{A}_q)<1$) can be readily proved based on \cite[Lemma~3]{Pang2020a}. For the second part of the result, the following lemma is invoked to facilitate the proof.
\begin{Lemma}
(see \cite[Lemma~5.6.10]{Horn1990}) Let $\rho(A)$ be the spectral radius of a (square) matrix $A$. For any given $\varrho > 0$, there exists a matrix norm $\|\cdot\|_E$ such that $\rho(A) \leq \|A\|_E \leq \rho(A)+\varrho$.
\end{Lemma}

Based on the above lemma, we can choose $\varrho \in (0,1-\max_{q\in\bar{\mathcal{V}}}\rho(\tilde{A}_q))$, then there exists a matrix norm $\|\cdot\|_E$ such that $\|\tilde{A}_q\|_E \leq \rho(\tilde{A}_q))+\varrho< 1, \forall q\in\bar{\mathcal{V}}$, which completes the proof.




\bibliographystyle{IEEEtran}
\bibliography{d_ne_coalition_grad_track_reference}

\end{document}